\documentclass[aos,MSNbibl,seceqn,rotating,citesort,dvips]{arximspdf}
\usepackage{subenv}

\doi{10.1214/11-AOS966}
\volume{40}
\issue{1}
\pubyear{2012}
\firstpage{436}
\lastpage{465}

\makeatletter

\newcommand{\xrightarrow}[1]{\stackrel{#1}{\rightarrow}}

\newtheorem{theorem}{Theorem}[section]

\newproclaim{Assump}{Assumption}

\newcommand{\sumtT}{\sum_{t=1}^T}
\newcommand{\sumiN}{\sum_{i=1}^N}

\newcommand{\vech}{\operatorname{vech}}
\newcommand{\veck}{\operatorname{veck}}

\renewcommand{\vec}{\operatorname{vec}}

\newcommand{\tr}{\operatorname{tr}}
\newcommand{\diag}{\operatorname{diag}}

\newtheorem{proposition}{Proposition}[section]

\makeatother

\begin{document}
\begin{frontmatter}

\title{Statistical analysis of factor models of~high~dimension}
\runtitle{Factor models of high dimension}

\begin{aug}
\author[A]{\fnms{Jushan} \snm{Bai}\corref{}\thanksref{t1}\ead[label=e1]{Jushan.bai@columbia.edu}}
\and
\author[B]{\fnms{Kunpeng} \snm{Li}\thanksref{t2}\ead[label=e3]{likp.07@sem.tsinghua.edu.cn}}
\runauthor{J. Bai and K. Li}
\affiliation{Columbia University and the Central University of Finance and Economics,
and Tsinghua University and~University~of International
Business and Economics}
\address[A]{Department of Economics\\
Columbia University\\
420 West 118th Street\\
New York, New York 10027\\
USA\\
\printead{e1}}
\address[B]{Department of Economics \\
School of Economics and Management \\
Tsinghua University \\
and\\
Department of Quantitative Economics \\
School of International Trade\\
\quad and Economics \\
University of International Business\\
\quad and Economics\\
Beijing, 100084\\
China\\
\printead{e3}} 
\end{aug}

\thankstext{t1}{Supported by NSF Grant SES-0962410.}

\thankstext{t2}{Supported by the new entrant program of UIBE
(Grant Number: 11QD11).}

\received{\smonth{6} \syear{2011}}
\revised{\smonth{10} \syear{2011}}

%
\begin{abstract}
This paper considers the maximum likelihood estimation of factor models
of high dimension, where the number of variables ($N$) is comparable
with or even greater than the number of observations~($T$). An
inferential theory is developed. We establish not only consistency but
also the rate of convergence and the limiting distributions. Five
different sets of identification conditions are considered. We show
that the distributions of the MLE estimators depend on the
identification restrictions. Unlike the principal components approach,
the maximum likelihood estimator explicitly allows
heteroskedasticities, which are jointly estimated with other
parameters. Efficiency of MLE relative to the principal components
method is also considered.
\end{abstract}

%

\end{frontmatter}

\section{Introduction}\label{sec1}

Factor models provide an
effective way of summarizing information from large data sets, and are
widely used in social and physical
sciences.\setcounter{footnote}{2}\footnote{See, for example,
\cite{r10,r12,r18,r19,r20,r21}.}
There has also been advancement in the theoretical analysis of factor
models of high dimension. Much of this progress has been focused on the
principal components method; see, for example,
\cite{r6,r8,r20} and~\cite{r21}.\footnote{\cite{r14} considers the principal components analysis in the frequency domain
for generalized dynamic factor models.} The advantage of the principal components
method is that it is easy to compute and it provides consistent
estimators for the factors and factor loadings when both $N$ and $T$
are large. The principal components method implicitly assumes that the
idiosyncratic covariance matrix is a~scalar multiple of an identity
matrix. While the method is robust to heteroscedasticity and weak
correlations in the idiosyncratic errors, there are biases associated
with the estimates. In fact, if $N$ is fixed, the principal components
estimator for the factor loadings is inconsistent, as shown in
\cite{r6}, except under homoscedasticity.

In this paper, we consider the maximum likelihood estimator under the
setting of large $N$ and large $T$. The maximum likelihood estimator
(MLE) is more efficient than the principal components method. In
addition, MLE is also consistent and efficient under fixed $N$ and
large $T$ because this setting falls within the framework of classical
inference; see, for example,~\cite{r4} and~\cite{r17}.

Our estimator coincides with the classical factor analysis. However,
the statistical theory does not follow from the existing literature.
Classical inferential theory is based on the assumption that $N$ is
fixed (or one of the dimensions is fixed). This assumption, to a~certain extent, runs counter to the primary purpose of factor analysis,
which is to explain the commonality among a~large number of variables
in terms of a~small number of latent factors.
Let $M_{z z}=\frac1 {T-1} \sum_{t=1}^T (z_t-\bar z)(z_t-\bar z)'$ be
the data matrix of $N\times N$ and let $\Sigma_{z z}(\theta)=E(M_{z
z})$. A~key assumption in classical inference is that
$\sqrt{T} \operatorname{vech} (M_{z z}-\Sigma_{z z}(\theta))$ is
asymptotically normal with a~positive definite limiting covariance
matrix, as $T\rightarrow\infty$. This assumption does not hold as $N$
also goes to infinity. The asymptotic normality is not well defined
with an increasing dimension. For example, if $N>T$, $M_{z z}$ is a~singular matrix, so it cannot have a~normal distribution with a~positive covariance matrix.
Furthermore, the dimension of the unknown parameters (denoted by
$\theta
$) is also increasing as
$N$ increases. The usual delta method (Taylor expansion) for deriving
the limiting distribution of the MLE of $\theta$ will not work.
Therefore, the high-dimensional inference for MLE requires a~new framework.

Fixing $N$ is for the purpose of tractability for theoretical
analysis. Such an assumption is unduly restrictive. Many applications
or theoretical models involve data sets with the number of variables
comparable with or even greater than the number of observations; see
\cite{r10,r19,r20} and~\cite{r21}. Although the
large-$N$ analysis is demanding, the limiting distribution of the
maximum likelihood estimator has a~much simpler form under large $N$
than under fixed~$N$.

There exists a~small literature on efficient estimation of factors and
factor loadings under large $N$.
\cite{r9} considers a~two-step approach by treating both the
factors and the factor loadings as the parameters of interest.
\cite{r11} also considers a~two-step
approach. The first step uses the principal components method to obtain
the residuals and the second step uses a~feasible generalized least
squares. This method depends on large $N$ and large $T$ to get
consistent estimation of the residual variances. MLE is considered by
\cite{r13}. 
A~certain average consistency is obtained;~\cite{r13}~does not
consider consistency for individual parameters nor the limiting distributions.

The present paper is the first to develop a~full statistical theory
for the maximum likelihood estimator. Our approach is different from
the existing literature.
The challenge of the analysis lies in the simultaneous estimation of
the heteroscedasticities and other model parameters. To estimate the
heteroscedasticity, the maximum likelihood estimator does not rely on
estimating the individual residuals, which would be the case for
two-step procedures. Using residuals to construct variance estimators
will be inconsistent when one of the dimension is fixed. The MLE
remains consistent under fixed $N$.

The rest of this paper is organized as follows. Section~\ref{sec2} introduces
the model and assumptions. Section~\ref{sec3} considers a~symmetrical
presentation for factor models. Identification conditions are
considered in Section~\ref{sec4}. Consistency and limiting distributions are
derived in Section~\ref{sec5}. Section~\ref{sec6} considers the estimation of factor
scores. Section~\ref{sec7} compares the efficiency of the MLE relative to the
principal components method and Section~\ref{sec8} discusses computational
issues. The last section concludes. Proofs of consistency are given in
the \hyperref[app]{Appendix} and additional proofs are provided in the
supplement
\cite{bai2011Supplement}. Throughout the paper, the norm of a~vector or
matrix is that of Frobenius, that is, $\|A\|=[\tr(A'A)]^{1/2}$ for
vector or matrix $A$; $\diag(A)$ represents a~diagonal matrix when $A$
is a~vector, but $\diag(A)$ can be either a~matrix or a~column vector
(consisting of the diagonal elements of $A$) when $A$ is a~matrix.

\section{Factor models}\label{sec2}

Let $N$ denote the number of variables and $T$ the sample size. For
$i=1,\ldots, N$ and $t=1,\ldots, T$, the
observation $z_{it}$ is said to have a~factor structure if it can be
represented as
%
%
\begin{equation}\label{sec21}
z_{it}=\alpha_{i}+\lambda_i'f_t + e_{it},
\end{equation}
where $f_t=(f_{t1},f_{t2},\ldots,f_{tr})'$ and $\lambda_i=(\lambda
_{i1},\ldots,\lambda_{ir})'$; both are $r\times1$. Let
$\Lambda=(\lambda_1,\lambda_2,\ldots,\lambda_N)'$ be $N\times r$, and
$z_t=(z_{1t},\ldots,z_{Nt})'$ be the $N\times1$ vector of observable
variables. Let $e_t$ and
$\alpha$ be similarly defined.
In matrix form,
%
%
\begin{equation}\label{sec22}
z_t=\alpha+ \Lambda f_t + e_t.
\end{equation}
The vector $z_t$ is observable; none of the right-hand side variables
are observable.
We make the following assumptions:
\renewcommand{\theAssump}{\Alph{Assump}}
%
\begin{Assump}\label{AssumptionA}
$\{f_t\}$ is a~sequence of fixed
constants. Let $M_{f f}=\frac1 T \times\sum_{t=1}^T(f_t-\bar
f)(f_t-\bar
f)'$ be the sample variance of $f_t$ where $\bar f=\frac{1}{T}\sum
_{t=1}^T f_t$.
There exists an $\overline M_{f f}>0$ (positive definite) such that
$\overline M_{f f}=\lim_{T \to\infty}M_{f f}. $
\end{Assump}
%
%
\begin{Assump}\label{AssumptionB}
$E(e_t)=0$; $E(e_te_t')=\Sigma_{ee}=\diag
(\sigma_1^{2},\sigma_2^{2},\ldots,\sigma_N^{2})$;\break \mbox{$E(e_{it}^4)\leq C^4$}
for all $i$ and $t$, for some $C<\infty$. The $e_{it}$ are independent
for all $i$ and $t$, and the $N\times1$ vector $e_t$ is identically
distributed over $t$.
\end{Assump}
%
%
\begin{Assump}\label{AssumptionC}
There exists a~positive constant $C$ large
enough such that:
\begin{longlist}[(C.3)]
\item[(C.1)] $\|\lambda_j\| \le C$ for all $j$.
\item[(C.2)] $C^{ - 2} \le\sigma_j^2 \le C^2$ for all $j$.
%
\item[(C.3)] The limits $\lim_{N \to\infty}
N^{-1}\Lambda'\Sigma_{ee}^{-1}\Lambda=Q$ and $ \lim_{N \to
\infty}
\frac{1}{N}\sum_{i=1}^N\sigma_i^{-4}(\lambda_i\otimes
\lambda_i)
(\lambda_i'\otimes\lambda_i')= \Omega$ exist, where $Q$ and
$\Omega$
are positive definite matrices.
\end{longlist}
\end{Assump}
%
%
\begin{Assump}\label{AssumptionD}
The variances $\sigma_j^2$ are estimated in
the compact set $[C^{-2}$, $C^2]$. Furthermore, $M_{f f}$ is restricted
to be in a~set consisting of all semi-positive definite matrices with
all elements bounded in the interval $[-C,C]$, where $C$ is a~large constant.
\end{Assump}

In Assumption~\ref{AssumptionA}, we assume $f_t$ is a~sequence of fixed constants. Our
analysis holds if it is a~sequence of random variables. In this case,
we assume $f_t$ to be independent of all other variables.
The analysis can then be regarded as conditioning on $\{f_t\}$. 
Without loss of generality, we assume that
$\bar f=\frac{1}{T}\sum_{t=1}^T f_t=0$ [or $E(f_t)=0$ for random
factors] because the model can be rewritten as
$z_t=\alpha+\Lambda\bar f+ \Lambda(f_t-\bar f) + e_t=\alpha
^*+\Lambda f_t^*+e_t$ with $\alpha^*=\alpha+\Lambda\bar f$ and
$f_t^*=f_t-\bar f$.
In Assumption~\ref{AssumptionB}, we assume $e_t$ to be independent over time. In fact,
our consistent result still holds if $e_t$ are serially correlated and
heteroscedastic over time or $e_{it}$ are correlated over $i$, provided
that these correlations are weak (sufficient conditions are given in
\cite{r3}). The limiting distribution then would need modification. For
simplicity, we shall consider the uncorrelated case. The analysis of
the maximum likelihood estimation under high dimension is already
difficult; allowing correlation will make the analysis even more
cumbersome. We will report the results under general correlation
patterns in a~separate paper.
Assumption~\ref{AssumptionD} is for theoretical analysis. Like all nonlinear
(nonconvex) analysis, parameters are assumed to be in a~bounded set.

The second moment of the sample, denoted by $M_{z z}$, is
%
%
\begin{equation}\label{sec23}
M_{z z} = \frac{1}{T}\sum_{t = 1}^T (z_t-\bar z)(z_t-\bar z)',
\end{equation}
where $\bar z=T^{-1}\sum_{t=1}^T{z_t}$. Note that the division by $T$
instead of $T-1$ is for notational simplicity.
Let $\Sigma_{z z}$ be
%
%
\begin{equation}\label{sec20}
\Sigma_{z z}=\Lambda M_{f f} \Lambda' +\Sigma_{ee}.
\end{equation}
The objective function considered in this paper is
%
%
\begin{equation}\label{sec26}
\ln L=- {\frac{1}{2N}\ln}|\Sigma_{z z}| - \frac{1}{2N}\tr( M_{z z}
\Sigma_{z z}^{ - 1}).
\end{equation}
The above objective function may be regarded as a~quasi likelihood
function. To see this, assume $f_t$ is stochastic with mean zero\vadjust{\goodbreak} and
variance $\Sigma_{f f}$. From $z_t=\alpha+\Lambda f_t+e_t$, the
variance matrix of $z_t$, denoted by $\Sigma_{z z}$, is
\[
\Sigma_{z z}=\Lambda\Sigma_{f f} \Lambda' +\Sigma_{ee}.
\]
So the quasi likelihood function (omitting a~constant) can be written as
\begin{eqnarray*}
\ln L &=&-{\frac{1}{2N}\ln}|\Sigma_{zz}|-\frac{1}{2NT}\sum_{t=1}^T(z_t
-\alpha)'\Sigma_{zz}^{-1}(z_t-\alpha)
\\
&=&- {\frac{1}{2N}\ln}|\Sigma_{zz}| - \frac{1}{2NT}\sum
_{t=1}^T(z_t-\bar
z)'\Sigma_{zz}^{-1}(z_t-\bar z)\\
&&{}- \frac{1}{2NT}\sum_{t=1}^T(\bar z -
\alpha)'\Sigma_{zz}^{-1}(\bar z - \alpha).
\end{eqnarray*}
Clearly $\hat\alpha$ minimizes the likelihood function at $\bar z$. So
the concentrated quasi likelihood function can now be written as
\begin{eqnarray*}
\ln L &=& - {\frac{1}{{2N}}\ln}|\Sigma_{z z}| - \frac{1}{{2N}}\tr
\Biggl[\frac{1}{T}\sum_{t = 1}^T (z_t - \bar z)(z_t - \bar z)'\Sigma_{z z}^{
- 1}\Biggr]\\
&=&- {\frac{1}{{2N}}\ln}|\Sigma_{z z}| - \frac{1}{2N}\tr( M_{z z}
\Sigma
_{z z}^{ - 1}),
\end{eqnarray*}
which is the same as~(\ref{sec26}) except that $\Sigma_{f f}$ is in
place of $M_{f f}$. Because the factors are fixed constants instead of
random variables, as stated in Assumption~\ref{AssumptionA}, it is natural to use
$M_{f f}$ rather than $\Sigma_{f f}$ in~(\ref{sec20}) and~(\ref{sec26}).

If both $\Lambda$ and $F=(f_1,f_2,\ldots,f_T)'$ are treated as parameters,
the corresponding likelihood function is
%
%
\begin{equation} \label{usual-mle}
-{\frac{1}{2N}\ln}|\Sigma
_{ee}|-\frac
{1}{2NT}\sumtT(z_t-\alpha-\Lambda f_t)'\Sigma_{ee}^{-1}(z_t-\alpha
-\Lambda f_t).
\end{equation}
%
Since $\Lambda$ has $Nr$ parameters and $F$ has $Tr$ parameters, the
number of parameters to be estimated will be very large, which leads to
efficiency loss. In contrast, the number of parameters in~(\ref{sec26})
is only $N(r+1)+r(r+1)/2$, which is considerably smaller than the
number of parameters in~(\ref{usual-mle}), which is $N(r+1)+Tr$. The
difference is pronounced for small $N$ but large $T$.
In fact, when estimating $\Lambda, F$ and $\Sigma_{ee}$, the global
maximum likelihood estimator does not exist. It can be shown that the
likelihood function diverges to infinity by certain choice of
parameters (see~\cite{r3}, page 587).

By restricting $\Sigma_{ee}=I_N$ (an identity matrix), the MLE
estimator of~(\ref{usual-mle}) becomes the principal components
estimator. That is, the principal components method minimizes the
objective function
$ \sumtT(z_t-\alpha-\Lambda f_t)'(z_t-\break\alpha-\Lambda f_t)$
over $\alpha$, $\Lambda$ and $F$.
The estimators cannot be efficient when heteroscedasticity actually exists.

Even though the $f_t$ are fixed constants, we avoid directly
estimating $f_t$. Instead we only estimate the sample moment of $f_t$.
This considerably reduces the number of parameters and removes the
corresponding incidental parameters bias. 
The estimator is also consistent under fixed $N$, since the setting
falls back to the classical factor analysis.

By maximizing~(\ref{sec26}), in combination with~(\ref{sec20}), we can
obtain three first-order conditions (see, e.g.,~\cite{r17}):
%
%
\begin{eqnarray}\label{sec271}
&\hat\Lambda' \hat\Sigma_{z z}^{ - 1}( M_{z z}-\hat\Sigma_{z z})=0,&
\\
%
%
\label{sec272}
&\diag( {\hat\Sigma_{z z}^{ - 1}} ) = \diag( \hat\Sigma_{z z}^{ -
1} M_{z z} \hat\Sigma_{z z}^{ - 1}),&
\\
%
%
\label{sec273}
&\hat\Lambda' \hat\Sigma_{z z}^{ - 1} \hat\Lambda=\hat\Lambda'
\hat\Sigma_{z z}^{ - 1} M_{z z} \hat\Sigma_{z z}^{ - 1} \hat
\Lambda,&
\end{eqnarray}
where $\hat\Lambda$, $\hat M_{f f}$ and $\hat\Sigma_{ee}$ denote the
MLE and $\hat\Sigma_{zz}=\hat\Lambda\hat M_{f f} \hat\Lambda'
+\hat\Sigma_{ee}$.

Condition~(\ref{sec271}) is derived from the partial derivatives with
respect to~$\Lambda$, (\ref{sec272}) is derived with respect to the
diagonal elements of $\Sigma_{ee}$, and~(\ref{sec273}) is derived with
respect to $M_{f f}$. Equation~(\ref{sec273}) can be obtained from
(\ref{sec271}) by post-multiplying $\hat\Sigma_{z z}^{ - 1} \hat
\Lambda$. Since~(\ref{sec273}) is redundant, in order to make the
system of three equations solvable, we need to impose further
restrictions. These identification restrictions
will be discussed in Section~\ref{sec4}.

\section{Symmetry and choice of representations}\label{sec3}
Consider the model
\[
z_{it}=\delta_t + \lambda_i'f_t + e_{it}.
\]
Let $z_i=(z_{i1},z_{i2},\ldots,z_{iT})'$, $\delta=(\delta_1,\ldots,\delta
_T)'$, $F=(f_1,f_2,\ldots,f_T)'$ and $e_i=(e_{i1},\ldots,e_{iT})'$; then
\[
z_i=\delta+F \lambda_i + e_i
\]
$ (i=1,2,\ldots,N)$.
Define
\[
M_{z z}=\frac1 N \sumiN(z_i-\bar z)(z_i-\bar z)',\qquad \Sigma_{z
z}=F M_{\lambda\lambda}F' + \Sigma_{ee}^\dag,
\]
where $\Sigma_{ee}^\dag=\diag(\sigma_1^2,\ldots,\sigma_T^2)$ and
\[
M_{\lambda\lambda}=\frac1 N \sumiN(\lambda_i-\bar\lambda
)(\lambda
_i-\bar\lambda)'
\]
is the $r\times r$ sample variance of the factor loadings.

Although we use the same notation of $M_{z z}$ and $\Sigma_{z z}$,
they are now
$T\times T$ matrices instead of $N\times N$. The matrix $\Sigma
_{ee}^\dag$ contains idiosyncratic variances in the time dimension
(time series heteroscedasticity). The quasi maximum likelihood estimator
maximizes the likelihood function
\[
\ln L= -{\frac1 {2T} \ln}|\Sigma_{z z}| -\frac1 {2T}\tr(M_{z
z}\Sigma
_{z z}^{-1}).
\]
%

This representation avoids estimating $\lambda_1,\lambda_2,\ldots
,\lambda_N$ directly, but only the sample moment of $\lambda_i$. The
representation has $Tr+T+r(r+1)/2$ number of parameters. If $N$ is much
larger than $T$, this representation will give more efficient
estimation. In particular, if $T$ is fixed, we can only use this
representation to get consistent estimation of $f_1,f_2,\ldots, f_T$
and $
\Sigma_{ee}^\dag$. This representation will also be useful if one is
interested in estimating the heteroscedasticity in the time dimension.

The analysis of one representation will carry over to the other by
switching the role of $N$ and $T$ and the role of $\Lambda$ and $F$. So
it is sufficient to carefully examine one representation.
Bearing this in mind, our analysis focuses on the representation in
the previous section.
The objective function in~(\ref{sec26}) involves fewer parameters when
$N$ is less than $T$, although we make no assumption about the relative
size between $N$ and $T$ (except for Theorem~\ref{thm7}), and in particular,
$N$ is allowed to be much larger than $T$.

\section{Identification conditions}\label{sec4}
It is well known that the factor models are not identifiable without
additional restrictions. For any\vspace*{1pt} $r\times r$ invertible matrix $R$, we
have $\Lambda M_{f f}\Lambda'= \tilde\Lambda\tilde M_{f f} \tilde
\Lambda'$ where $\tilde\Lambda=\Lambda R$ and $\tilde M_{f
f}=R^{-1}M_{f f} R^{'-1}$. Thus observationally equivalent models are obtained.
In order to uniquely fix $\Lambda$ and $M_{f f}$ given $\Lambda M_{f
f}\Lambda'$, we need
$r^2$ restrictions since an invertible $r\times r$ matrix has $r^2$
free parameters.
For details of identification conditions, readers are referred to
\cite{r17} and~\cite{r5}. There are many ways to impose restrictions. In
this paper, we consider five identification strategies which have been
used in traditional factor analysis.
These restrictions are listed in Table~\ref{table1}.
The left pane is for the representation in Section~\ref{sec2}, while the right
pane is for the representation in Section~\ref{sec3}.

%
%
\begin{sidewaystable}
\textwidth=\textheight
\tablewidth=\textwidth
\caption{Identifying restrictions}\label{table1}
\begin{tabular*}{\tablewidth}{@{\extracolsep{\fill}}lccccc@{}}
\hline
& \textbf{Restrictions on} $\bolds{F}$ & \textbf{Restrictions on}
$\bolds{\Lambda}$ & & \textbf{Restrictions on} $\bolds{\Lambda}$
& \textbf{Restrictions on} $\bolds{F}$\\
\hline
$\mbox{IC1}$ & Unrestricted & $\Lambda=(I_r,\Lambda_2')'$ & $\mbox{IC1}'$
& Unrestricted & $F=(I_r,F_2')'$\\
[3pt]
$\mbox{IC2}$
& $M_{f f} = \mbox{diagonal}$ & $\frac{1}{N}\Lambda'\Sigma_{ee}^{ -
1}\Lambda=I_r $&$\mbox{IC2}'$
& $M_{\lambda\lambda}=\mbox{diagonal}$
& $\frac{1}{T}F'\Sigma_{ee}^{\dag-1}F=I_r$\\
& (with distinct elements) &&& (with distinct elements)\\
[3pt]
$\mbox{IC3}$ & $ M_{f f} = I_r $
& $\frac{1}{N}\Lambda'\Sigma_{ee}^{ - 1}\Lambda=\mbox{diagonal}$
& $\mbox{IC3}'$ & $M_{\lambda\lambda}=I_r$
& $\frac{1}{T}F'\Sigma_{ee}^{\dag-1}F=\mbox{diagonal}$\\
&& (with distinct elements) &&& (with distinct elements)\\
[5pt]
$\mbox{IC4}$ & $M_{f f} = \mbox{diagonal}$
& $\Lambda=(\Lambda_1',\Lambda_2')'$
& $\mbox{IC4}'$ & $ M_{\lambda\lambda}= \mbox{diagonal}$
& $F=(F_1',F_2')'$\\[3pt]
&& $\Lambda_1
=\pmatrix{ 1 & 0 & \cdots& 0\cr
\lambda_{21} & 1 & \cdots& 0\cr
\vdots& \vdots& \ddots&\vdots\cr
\lambda_{r1}& \lambda_{r2}& \cdots& 1}$
&&& $F_1 =\pmatrix{ 1 & 0 & \cdots& 0 \cr
f_{21} & 1 & \cdots& 0\cr
\vdots& \vdots& \ddots&\vdots\cr
f_{r1}& f_{r2}& \cdots& 1}$
\\\\[-4pt]
$\mbox{IC5}$ & $ M_{f f} = I_r $ & $\Lambda'
=(\Lambda_1', \Lambda_2')'$
& $\mbox{IC5}'$ & $M_{\lambda\lambda}=I_r$ &
$F' =(F_1',F_2')'$\\[3pt]
&& $\Lambda_1 =\pmatrix{ \lambda_{11} & 0 & \cdots& 0 \cr
\lambda_{21} & \lambda_{22} & \cdots& 0 \cr
\vdots& \vdots& \ddots& \vdots\cr
\lambda_{r1} & \lambda_{r2} & \cdots& \lambda_{rr}}$
&&& $F_1 =\pmatrix{f_{11} & 0 & \cdots& 0 \cr
f_{21} & f_{22} & \cdots& 0 \cr
\vdots& \vdots& \ddots&\vdots\cr
f_{r1} & f_{r2} & \cdots& f_{rr}}$
\\\\[-6pt]
&& $\lambda_{ii}\neq0$, $i=1,2,\ldots, r$
&&& $f_{ii}\neq0$, $i=1,2,\ldots, r$\\
\hline
\end{tabular*}
\end{sidewaystable}

We make some comments on these restrictions. Given $\Lambda M_{f
f}\Lambda'$, IC1 will uniquely fix $\Lambda$ and $M_{f f}$. So full
identification is achieved.
But this is not the case for IC2. If we change the sign of any column
of $\Lambda$, $\Lambda M_{f f} \Lambda'$ is not changed. This implies
that we only identify $\Lambda$ up to a~column sign change.

Furthermore, if we switch the positions between the $i$th and $j$th
columns of $\Lambda$, and the positions between the $i$th and $j$th
diagonal elements of $M_{f f}$, the matrix $\Lambda M_{f f} \Lambda'$
is not changed. This means that we need restrictions on the ordering
of the diagonal elements of $M_{f f}$. In this paper, we assume that
the diagonal components of $M_{f f}$ are arranged from the largest to
the smallest and they must be distinct and positive. Because of this
restriction,
we naturally require that the diagonal elements of estimator
$\hat M_{f f}$ are also arranged in this order, which is important for
the proof of consistency.

Under IC3, for the same reason, we assume that the diagonal elements of
$\frac{1}{N}\Lambda'\Sigma_{ee}^{-1}\Lambda$ are distinct and positive,
and are arranged in decreasing order;
$\Lambda$~is identified up to a~column sign change.

IC4 imposes $\frac{1}{2}r(r+1)$ restrictions on the factor loadings,
and $\frac{1}{2}r(r-1)$ restrictions on the factors. Identification is
fully achieved like IC1.

Under IC5, we can only identify $\Lambda$ up to a~column sign change.
In addition,
we need nonzero diagonal elements for the lower triangular matrix. The
reason is intuitive. If the $i$th diagonal element is zero, both the
$i$th and $(i+1)$th columns will share the same structure.

IC1 is related to the measurement error problem; it assumes that the
first $r$ observations are noise measurements of the underlying
factors. IC2 and IC3 are the usual restrictions for MLE; see
\cite{r17}. IC4 and IC5 assume a~recursive relation: the first factor
affects the first variable only, and the first two factors affect the
first two variables only, and so on; they are widely used, for example,
\cite{r5} and~\cite{gewekezhou96}.
Clearly, IC1, IC4 and IC5 require a~careful choice of the first $r$
observations in practice. The inferential theory assumes that the
underlying parameters satisfy the restrictions, implying different
$\lambda_i$ under different restrictions.

\section{Asymptotic properties of the likelihood estimators}\label{sec5}

Since the number of parameters increases as $N$ increases, the usual
argument that the objective function converges in probability to a~fixed nonrandom function and the function achieves its maximum value at
the true parameter values will not work. This is because as $N$ and $T$
increase, there will be an infinite number of parameters in the limit.
Our idea of consistency is to obtain some average consistency, and then
use these initial results to obtain consistency for individual
parameters. Even the average consistency requires a~novel argument in
the presence of an increasing number of parameters.
%
%
\begin{proposition}\label{thm1}
Let $\hat\theta$ be the MLE by maximizing~(\ref{sec26}), where $\hat
\theta=(\hat\lambda_1,\ldots,\hat\lambda_N,\hat\sigma
_1^2,\ldots,\hat
\sigma_N^2,\hat M_{f f})$. Under Assumptions~\ref{AssumptionA}--\ref{AssumptionD}, when
$N,T\rightarrow
\infty$, with any one of the identification conditions \textup{IC1--IC5},
we have
%
%
\begin{subequation}\label{sec41}
%
%
\begin{eqnarray}
\label{sec411}
\frac{1}{N}\sum_{i=1}^N{\frac{1}{\hat\sigma_i^2}\|\hat\lambda
_i-\lambda
_i\|^2}&\xrightarrow{p}&0,\\
\label{sec412}
\frac{1}{N}\sum_{i=1}^N(\hat\sigma_i^2-\sigma_i^2
)^2&\xrightarrow{p}&0,\\
\label{sec413}
\hat M_{f f}- M_{f f}&\xrightarrow{p}&0.
\end{eqnarray}
\end{subequation}
\end{proposition}

Establishing the above result requires a~considerable amount of work.
Developing and identifying appropriate strategies have taken an even
greater amount of efforts. The difficulty lies in the problem of
infinite number of parameters in the limit and the nonlinearity of
objective function. The infinite number of parameters problem in this
paper is fundamentally different from those in the existing literature.
For example, consider an $AR(\infty)$ process $X_t=\sum_{j=1}^\infty
a_j\varepsilon_{t-j}$. Although there exist an infinite number of
parameters $\{a_j\}_{j=1}^\infty$, the assumption that $a_j\to0$, as
$j\to\infty$, effectively limits the number of parameters. For example,
$\hat a_j\equiv0$ is consistent for $a_j$ for $j\ge\ln(T)$. The
assumption that $a_j\rightarrow0$ may be viewed as one form of
smoothing restriction.
However, in the present context and in the absence of any form of
smoothness, all parameters are free parameters, and there will be an
infinite number of them in the limit.
This is the source of difficulty.

While there is also an infinite number of parameters problem in the
analysis of the principal components (PC) estimator, the method does
not estimate heteroscedasticity, and it minimizes an objective function
stated in Section~\ref{sec2}.
Its degree of nonlinearity is much less than the likelihood function~(\ref{sec26}).
It is the joint estimation of heteroscedasticity that makes the
analysis difficult.
In the \hyperref[app]{Appendix}, we provide a~novel proof of
consistency, which
constitutes a~departure from the usual analysis, say, in
\cite{NeMcF94} and~\cite{Vaart98}.


The proofs\vspace*{1pt} of~(\ref{sec411}) and~(\ref{sec413}) depend heavily on
the identification conditions.
If we denote $A\equiv(\hat\Lambda-\Lambda)'\hat\Sigma
_{ee}^{-1}\hat
\Lambda(\hat\Lambda'\hat\Sigma_{ee}^{-1}\hat\Lambda)^{-1}$, the proof
of consistency centers on proving $A\xrightarrow{p}0$. However, the
proof of $A\xrightarrow{p}0$ is quite different with different
identification conditions. Under IC2, for example, $(\hat\Lambda'\hat
\Sigma_{ee}^{-1}\hat\Lambda/N)^{-1}=I_r$. Under other identification
conditions, the proof of even $(\hat\Lambda'\hat\Sigma
_{ee}^{-1}\hat
\Lambda/N)^{-1}=O_p(1)$ is extremely demanding.
Under IC2, IC3 and IC5, we need to assume that the estimator $\hat
\Lambda$ has the same column signs as those of $\Lambda$ in order to
have consistency. Having the same column signs is regarded as part of
the identification restrictions under IC2, IC3 and IC5.

In order to derive the inferential theory for the estimated
parameters, we need to strengthen Proposition~\ref{thm1}. We state the
result as a~theorem:
%
%
\begin{theorem}\label{thm2}
Under the assumptions of Proposition~\ref{thm1}, we have
%
%
\begin{subequation}\label{sec42}
%
%
\begin{eqnarray}
\label{sec421}
\frac{1}{N}\sum_{i=1}^N\frac{1}{\hat\sigma_i^2}\|\hat\lambda
_i-\lambda
_i\|^2&=&O_p(T^{-1}),\\
\label{sec422}
\frac{1}{N}\sum_{i=1}^N(\hat\sigma_i^2-\sigma_i^2
)^2&=&O_p(T^{-1}),\\
\label{sec423}
\|\hat M_{f f}- M_{f f}\|^2&=&O_p({T^{-1}}).
\end{eqnarray}
\end{subequation}
\end{theorem}

It is interesting to compare our results with those in classical factor
analysis. If $N$ is fixed, the existing literature\vspace*{-1pt} has already shown
that $\hat\lambda_j$ and $\hat\sigma_j^2$ converge to $\lambda_j$ and
$\sigma_j^2$ at the rate of $\sqrt{T}$ for any $j$. Since\vspace*{-1pt} $N$ is fixed,
the classical result implies~(\ref{sec421}) and~(\ref{sec422}). In
fact $\|\hat M_{f f}- M_{f f}\|^2=O_p({T^{-1}})$ holds also since it
can be derived from the first two (the results analogous to~(\ref{sec423}) under IC1 when $N$ is finite can be seen in~\cite{r2}).
Theorem~\ref{thm2} shows that these results still hold in the large-$N$
setting despite estimating an increasing number of elements. However,
we point out that the rate stated in Theorem~\ref{thm2} is not the
sharpest. If IC2 or IC3 is adopted as identification conditions, then
$\|\hat M_{f f}- M_{f f}\|^2=O_p({T^{-1}})$ can be refined as $\|\hat
M_{f f}- M_{f f}\|^2=O_p({N^{-1}T^{-1}})+O_p(T^{-2})$. Because~(\ref{sec423}) is sufficient for the inferential theory to be developed, we
only state this general result.

As pointed out earlier, the behavior of $A\equiv(\hat\Lambda-\Lambda
)'\hat\Sigma_{ee}^{-1}\hat\Lambda(\hat\Lambda'\hat\Sigma
_{ee}^{-1}\hat
\Lambda)^{-1}$ is important in establishing consistency. In fact,
matrix $A$ plays a~key role in the inferential theory as well. The
convergence rate of $A$ depends on identification conditions. We use
$A_k$ in place of $A$ under IC$k$ ($k=1,2,\ldots,5$). Under IC2 and
IC3, the convergence rate of $A$ is $\min(\sqrt{NT},T)$. However, under
other sets of identification conditions, the convergence rate of $A$ is
$\sqrt{T}$. This difference in convergence rate affects the limiting
distribution of~$\hat M_{f f}$, which also makes the limiting
distributions of~$\hat\lambda_j$ different.

In Section C of the supplement~\cite{bai2011Supplement}, we give the asymptotic representations
of $\sqrt T(\hat\lambda_j-\lambda_j)$ under IC1--IC5. The main
representations are given in (C.5), (C.12), (C.17) and (C.24),
respectively. The following theorem is a~consequence of these representations.
%
%
\begin{theorem}\label{thm3}
Under the assumptions of Proposition~\ref{thm1}, for each
$j=1, 2,\ldots, N$, as $N, T\to\infty$, we have:

Under \textup{IC1},
%
%
\begin{subequation}\label{sec43}
%
%
\begin{eqnarray}
\label{sec431}
\sqrt{T}(\hat\lambda_j-\lambda_j)&\xrightarrow{d}&\mathcal
{N}\bigl(0,(\overline M_{f f})^{-1}(\lambda'_j\Sigma_{eer}\lambda_j+\sigma
_j^2)\bigr).
\end{eqnarray}

Under \textup{IC2} or \textup{IC3},
%
\begin{eqnarray}\label{sec432}
\sqrt{T}(\hat\lambda_j-\lambda_j)&\xrightarrow{d}&
\mathcal{N}(0,(\overline M_{f f})^{-1}\sigma_j^2).
\end{eqnarray}

Under \textup{IC4}, for $j>r$
%
\begin{eqnarray}\label{sec433}
\sqrt T(\hat\lambda_j-\lambda_j)
&\xrightarrow{d}&\mathcal
{N}\bigl( 0,\Pi+(\overline M_{f f})^{-1}\sigma_j^2\bigr)
\end{eqnarray}
and for $2\le j\le r$
\begin{eqnarray}
\sqrt T(\hat\lambda_j-\lambda_j)&\xrightarrow{d}&\mathcal{N}
\bigl(0,\Pi+2I_r^{j-1}(\overline M_{f f})^{-1}\sigma_j^2-(\overline M_{f
f})^{-1}\sigma_j^2\bigr).\nonumber
\end{eqnarray}

Under \textup{IC5}, for $j>r$
%
\begin{eqnarray}\label{sec434}
\sqrt T(\hat\lambda_j-\lambda
_j)&\xrightarrow
{d}&\mathcal{N}(0,\Xi+ I_r \sigma_j^2)
\end{eqnarray}
and for $1\le j\le r$
\begin{eqnarray}
\sqrt T(\hat\lambda_j-\lambda_j)&\xrightarrow
{d}&\mathcal
{N}(0,\Xi+2I_r^{j}\sigma_j^2- I_r \sigma_j^2),\nonumber
\end{eqnarray}
\end{subequation}
where\vspace*{-1pt} $\Pi=(\lambda_j'\otimes I_r)\tilde D(\overline M_{f f} )\Phi
\tilde D(\overline M_{f f})'(\lambda_j \otimes I_r)$, $\Xi=(\lambda'
_j \otimes I_r)\overline D \Gamma\overline D'(\lambda_j \otimes
I_r)$. $\Sigma_{eer}$ is an $r\times r$ diagonal matrix with the $j$th
diagonal element $\sigma_j^2$; $I_r^j$ is an $r\times r$ diagonal
matrix with the first $j$ diagonal elements being 1 and the rest being
0. The meanings of $\tilde D(\overline M_{f f}), \overline D, \Phi$
and $\Gamma$ are explained below.
\end{theorem}

The matrix $\tilde D(M)$ (with $M=\overline M_{f f}$) is a~generalized
duplication matrix of $r^2 \times\frac{1}{2}r(r+1)$ depending on the
diagonal matrix $M$; $\tilde D(M)$ can be constructed row by row in the
following way. Given the number $k$, $1\le k\le r^2$, we denote $j
=\lfloor{(k - 1)/r}\rfloor+1 $ and $i = k -(j-1)r$, where $\lfloor
\cdot\rfloor$ denotes the largest integer no greater than the
argument. If $i \ge j$, all elements of the $k$th row are zero, except
that the $(\frac{1}{2}(2r-j+2)(j-1)+i-j+1)$th element is 1; if $i<j$,
all elements are zero, except that the $(\frac
{1}{2}(2r-i+2)(i-1)-i+j+1)$th element is $-m_jm_i^{-1}$, where $m_j$ is
the $j$th diagonal element of $M$. The $r^2\times\frac{1}{2}r(r-1)$
matrix $\overline{D}$ under IC5 is also a~generalized duplication
matrix. Let $A$ be a~skew-symmetric matrix and let $\veck(A)$ be the
operator\vspace*{1pt} that stacks the elements of $A$ strictly below the diagonal
into a~vector (excluding diagonal elements). Then $\overline D$ is
defined as $\vec(A)=\overline D \veck(A)$.

Here are some examples for $\tilde D(M)$. If $M$ is a~scalar, then
$\tilde D(M)=1$. If $M=\diag(m_1,m_2)$, then
\[
\tilde D(M)=\pmatrix{1&0&0\cr
0&1&0\cr 0&-m_2/m_1&0\cr 0&0&1}.
\]
If $M=\diag(m_1,m_2,m_3)$, then
\[
\tilde D(M)=\pmatrix{
1&0&0&0&0&0\cr 0&1&0&0&0&0\cr 0&0&1&0&0&0\cr
0&-m_2/m_1&0&0&0&0\cr
0&0&0&1&0&0\cr 0&0&0&0&1&0\cr 0&0&-m_3/m_1&0&0&0\cr
0&0&0&0&-m_3/m_2&0\cr
0&0&0&0&0&1}.
\]

Here are some examples for $\overline D$, which only depends on the
dimension of $\overline M_{f f}$ (i.e., the number of factors). If
$r=1$, then $\overline D=0$.
If $r=2$, then $\overline D=(0,1,-1,0)'$. If $r=3$, then
\[
\overline D=\pmatrix{
0&0&0\cr 1&0&0\cr 0&1&0\cr
-1&0&0\cr 0&0&0\cr 0&0&1\cr 0&-1&0\cr 0&0&-1\cr 0&0&0}.
\]

The matrix $\Phi$ in~(\ref{sec433}) is the limiting covariance of $\vech
(A_4)$, where $A_4$ is the $A$ matrix defined earlier under IC4. The
asymptotic representation of $A_4$ is given by (C.15) in Section C of
the supplement~\cite{bai2011Supplement}. From this asymptotic representation, the elements of
$\Phi$ can be easily computed and are given by~(C.16). The matrix
$\Gamma$ in~(\ref{sec434}) is the limiting covariance of~$\veck(A_5)$,
where $A_5$ is the matrix $A$ under IC5, and is asymptotically
skew-symmetric. The asymptotic representation of $A_5$ is given by
(C.22). The elements of $\Gamma$ are determined by (C.23) in Section C
of the supplement~\cite{bai2011Supplement}.

\textit{Remarks}: In classical factor analysis, the MLE usually imposes
the restriction of IC3. The limiting distribution $\hat\lambda_j$ in
classical factor analysis (fixed~$N$) is very complicated; see
\cite{r17}. Most textbooks on multivariate statistics do not even present
the limiting distributions owing to its complexity. As pointed out by
Anderson (\cite{r3}, page 583), the limiting distribution is ``too
complicated to derive or even present here.'' In contrast, the limiting
distribution under IC3 with large $N$ is
\[
\sqrt{T}(\hat\lambda_j-\lambda_j) \xrightarrow{d}\mathcal
{N}(0,(\overline M_{f f})^{-1}\sigma_j^2).
\]
This is as efficient as the case in which the $f_t$ are observable. For
if the~$f_t$ are observable, the estimator of $\lambda_j$ by applying
OLS to~(\ref{sec21}) is $\hat\lambda_j^{\mathrm{ols}}=(T^{-1}\sum_{t=1}^T
(f_t-\bar f)(f_t-\bar f)')^{-1}(T^{-1}\sum_{t=1}^T(f_t-\bar
f)(z_{jt}-\bar z_j))$.
It is easy to show that $\sqrt{T}(\hat\lambda_j^{\mathrm{ols}}-\lambda
_j)\xrightarrow{d}\mathcal{N}(0,(\overline M_{f f})^{-1}\sigma_j^2)$,
the same as the MLE estimator under IC2 and IC3. The OLS estimator is
the best linear unbiased estimator under Assumption~\ref{AssumptionB}, as the error
terms of the regression equation are i.i.d.


Now we state\vspace*{-1pt} the limiting distributions of $\hat M_{f f}$ and $\hat
\sigma_j^2$ for each $j$. In Section C of the supplement~\cite{bai2011Supplement}, we derive the
asymptotic representations for $\hat M_{f f}-M_{f f}$ under IC1, IC2
and IC4, respectively, which are given by~(C.7), (C.11) and (C.19). The
asymptotic representation for $\hat\sigma_j^2-\sigma_j^2$ is given by~(C.4).
The next two theorems follow these asymptotic representations.
%
%
\begin{theorem}\label{thm4}
Under the assumptions of Proposition~\ref{thm1}, we have:
%
%

Under \textup{IC1},
\[
\sqrt{T} \vech(\hat M_{f f}-M_{f f})
\xrightarrow{d}\mathcal{N}\bigl(0,4D_r^+(\Sigma_{eer} \otimes\overline M_{f
f})D_r^{+\prime}\bigr).
\]

Under \textup{IC2}, $N/T\rightarrow0$ and normality of $e_{it}$,
\begin{eqnarray*}
&&\sqrt{NT}\diag(\hat M_{f f}-M_{f f})\\
&&\qquad\xrightarrow{d}
\mathcal{N}\bigl(0,J_r[2(I_r\otimes\overline M_{f f})\Omega(I_r\otimes
\overline M_{f f})+4(Q\otimes\overline M_{f f})]J_r'\bigr).
\end{eqnarray*}

Under \textup{IC4},
\[
\sqrt{T}\diag(\hat M_{f f}-M_{f f})
\xrightarrow{d}\mathcal{N}\bigl(0,4J_r[(\Lambda_1'\Sigma_{eer}^{-1}\Lambda
_1)^{-1}\otimes\overline M_{f f}]J_r'\bigr),
\]
where $D_r^{+}$ is the Moore--Penrose inverse of the duplication matrix
$D_r$; $J_r$ is an $r\times r^2$ matrix, which satisfies, for any
$r\times r$ matrix M, $\diag\{M\}=J_r \vec(M)$, where $\diag\{\cdot
\}$
is the operator which stacks the diagonal elements into a~vector.
\end{theorem}

Note under IC3 and IC5, $M_{f f}$ is known and thus not estimated.
Normality under IC2 is used only for calculating the limiting variance.
Given the asymptotic representation of $\hat M_{f f}-M_{f f}$, it is
easy to derive the limiting distribution under nonnormality.
%
%
\begin{theorem}\label{thm5}
Under the assumptions of Proposition~\ref{thm1}, with any set of the
identification conditions, we have
%
%
\begin{equation}\label{sec48}
\sqrt{T}(\hat\sigma_j^{2}-\sigma_j^{2})\xrightarrow{d}\mathcal
{N}\bigl(0,\sigma_j^{4}(2+\kappa_j)\bigr),
\end{equation}
where $\kappa_j$ is the excess kurtosis of $e_{jt}$. Under normality of
$e_{it}$, the limiting distribution becomes $\mathcal{N}(0,2\sigma_j^4)$.
\end{theorem}

Our analysis assumes that the underlying parameters satisfy the
identification restrictions,
which is also the classical framework of~\cite{r17}. 
A~consequence is that we are directly estimating the underlying true
parameters instead of rotations of them. The rotation matrix used in
\cite{r21} and~\cite{r8} degenerates into an identity matrix. This
result itself is interesting.

\section{Asymptotic properties for the estimated factors}\label{sec6}
The factors $f_t$ can be estimated by two different methods. One is the
projection formula and the other is the generalized least squares
(GLS). These methods are discussed in~\cite{r3}.

If the factor $f_t$ is normally distributed with mean zero and variance
$\Sigma_{f f}$, and is independent of $e_t$, then the joint
distribution of $(f_t,z_t)$, by~(\ref{sec22}), can be written as
%
%
\begin{equation}\label{sec51}
\left[\matrix{
f_t \cr
z_t \cr}\right]
\sim N \left[\pmatrix{0\cr\alpha},\matrix{
\Sigma_{f f}&\Sigma_{f f}\Lambda'\cr
\Lambda\Sigma_{f f}&\Lambda\Sigma_{f f}\Lambda'+\Sigma_{ee}}
\right].
\end{equation}
Given $z_t$, the best predictor of $f_t$, $f_t^{p}$, is $f_t^{p} =
\Sigma_{f f}\Lambda'(\Lambda\Sigma_{f f}\Lambda'+\Sigma
_{ee})^{-1}(z_t-\alpha)$. By the basic result $(\Lambda\Sigma_{f
f}\Lambda'+\Sigma_{ee})^{-1}=\Sigma_{ee}^{-1}-\Sigma
_{ee}^{-1}\Lambda
(\Sigma_{f f}^{-1}+\Lambda'\Sigma_{f f}^{-1}\Lambda)^{-1}
\times\Lambda
'\Sigma_{ee}^{-1}$, we have
%
%
\begin{equation}\label{sec52}
f_t^{p}=(\Sigma_{f f}^{-1}+\Lambda'\Sigma_{ee}^{-1}\Lambda
)^{-1}\Lambda
'\Sigma_{ee}^{-1}({z_t-\alpha}).
\end{equation}
Although~(\ref{sec52}) is deduced under the assumption of normality of
$f_t$ and in this paper the $f_t$ are fixed constants, equation~(\ref{sec52}) can still be used to estimate $f_t$ by replacing the
parameters with their corresponding estimates. So the estimator $\tilde
f_t$ is
%
%
\begin{equation}\label{sec53}
\tilde f_t=(\hat M_{f f}^{-1}+\hat\Lambda'\hat\Sigma_{ee}^{-1}\hat
\Lambda)^{-1}\hat\Lambda'\hat\Sigma_{ee}^{-1}({z_t-\bar z}).
\end{equation}

An alternative procedure is the GLS. If the $\lambda_j$ and $\sigma
_j^2$ are observable, the GLS estimator of $f_t$ is $(\Lambda'\Sigma
_{ee}^{-1}\Lambda)^{-1}\Lambda'\Sigma_{ee}^{-1}({z_t-\bar z})$. The
unknown variables can be replaced by their estimates. We define the GLS
estimator of $f_t$ as
%
%
\begin{equation}\label{sec54}
\hat f_t=(\hat\Lambda'\hat\Sigma_{ee}^{-1}\hat\Lambda)^{-1}\hat
\Lambda
'\hat\Sigma_{ee}^{-1}({z_t-\bar z}).
\end{equation}
Under large $N$, not much difference exists between
(\ref{sec53}) and~(\ref{sec54}). In fact, they are asymptotically
equivalent and have the same limiting distributions.
But for relatively small $N$, the difference may not be ignorable.
%
%
\begin{proposition}\label{thm6}
Under the assumptions of Proposition~\ref{thm1}, $\tilde f_t=\hat
f_t+O_p(1/N)$.
\end{proposition}

Since $\tilde f_t$ and $\hat f_t$ have the same limiting distribution,
we only state the distribution for~(\ref{sec54}). In Section D of the
supplement~\cite{bai2011Supplement} we derive the asymptotic representations for $\sqrt N(\hat
f_t-f_t)$, which are given by (D.3), (D.4) and~(D.5), respectively.
From these representations, we obtain:
%
%
\begin{theorem}\label{thm7}
Let $\Delta\in[0,\infty)$. Under the assumptions of Proposition~\ref{thm1} and $\sqrt N/T\to0$, we have:

Under \textup{IC1} and $N/T\rightarrow\Delta$,
%
%
\begin{subequation}\label{sec55}
%
\begin{eqnarray}
\label{sec551}
\sqrt N(\hat f_t -f_t)\xrightarrow{d}\mathcal{N}
\bigl(0,\Delta f_t'(\overline M_{f f})^{-1}f_t\Sigma_{eer}+Q^{-1}\bigr).
\end{eqnarray}

Under \textup{IC2},
%
\begin{eqnarray}\label{sec552}
\sqrt{N}(\hat f_t-f_t)&\xrightarrow{d}&\mathcal{N}(0,I_r).
\end{eqnarray}

Under \textup{IC3},
%
\begin{eqnarray}\label{sec553}
\sqrt{N}(\hat f_t-f_t)&\xrightarrow{d}&
\mathcal{N}(0,Q^{-1}).
\end{eqnarray}

Under \textup{IC4} and $N/T\rightarrow\Delta$,
%
\begin{eqnarray}\label{sec554}
&&\sqrt N(\hat f_t-f_t)\nonumber\\[-8pt]\\[-8pt]
&&\qquad\xrightarrow{d}\mathcal{N}\bigl(0,\Delta(I_r
\otimes f_t')\tilde D(\overline M_{f f})\Phi\tilde D'(\overline M_{f
f})(I_r\otimes f_t)+Q^{-1}\bigr).\nonumber
\end{eqnarray}
Under \textup{IC5} and $N/T\rightarrow\Delta$,
%
\begin{eqnarray}\label{sec555}
\sqrt{N}(\hat f_t-f_t)&\xrightarrow{d}& \mathcal{N}
\bigl(0,\Delta(I_r\otimes f_t')\overline D\Gamma\overline D'(I_r\otimes
f_t) + Q^{-1}\bigr).
\end{eqnarray}
\end{subequation}
The matrix $Q$ is defined in Assumption~\ref{AssumptionC}. The matrices $\Sigma_{eer}$,
$\Phi,\Gamma,\tilde D(M_{f f})$ and $\overline D$ are defined in
Theorem~\ref{thm3}.
\end{theorem}

If $\Delta=0$, Theorem~\ref{thm7} shows that $\hat f_t$ has the same
limiting distribution regardless of the identification restrictions. In
this case, the variance is equal to $Q=\lim_{N\to
\infty
}N^{-1}\Lambda'\Sigma_{ee}^{-1}\Lambda$. Note that under IC2, $Q=I_r$.
Irrespective of whether $\Delta$ is zero,\vadjust{\goodbreak}
$\hat f_t$ is efficient under IC2 and IC3 in the sense that the
limiting variance coincides with the situation in which both the factor
loadings and the variances $\sigma_t^2$ are observable and GLS is
applied to a~cross-sectional regression for each fixed $t$. This result
requires $\sqrt{N}/T\rightarrow0$.

Recall that a~factor model has two symmetrical representations. We
also discussed earlier which presentation should be used in practice.
If $N$ is smaller than~$T$, we should estimate the factor loadings by
the maximum likelihood method because this representation has fewer
number of parameters. The opposite is true if $T$ is smaller than $N$.
This intuitive argument is borne out by the results of Theorems~\ref{thm3} and~\ref{thm7}.
By\vspace*{1pt} Theorem~\ref{thm3}, the magnitudes of $N$ and $T$ do not affect the
limiting covariance of $\hat\lambda_j$ (other than the rate of
convergence) but they do affect the limiting covariance of $\hat f_t$
as shown by Theorem~\ref{thm7}. If $\Delta$ is large, $\hat f_t$ cannot
be estimated well under IC1, IC4 and IC5. Note that $\hat f_t$ is not
the maximum likelihood estimator. In this case, we can use the
representation of Section~\ref{sec3} and directly estimate~$f_t$ by the maximum
likelihood method.

\section{Comparison with the principal components method}\label{sec7}

The method of principal components (PC) does not assume a~factor model,
and the method is usually regarded as a~dimension reduction technique.
But PC can be used to estimate factor models under large $N$ and large
$T$; see~\cite{r12,r10,r8} and~\cite{r20}. Let\vspace*{-1pt}
$\hat
\lambda_j^{\mathrm{pc}}$ and $\hat f_t^{\mathrm{pc}}$ denote the PC estimators for
$\lambda_j$ and $f_t$, respectively. The results of~\cite{r6} and
\cite{bai-ng-2010} imply the following asymptotic representation for the
principal components estimators. As $N,T\rightarrow\infty$, if $\sqrt
N/T \to0$, then
\begin{eqnarray*}
\sqrt T(\hat\lambda_j^{\mathrm{pc}}-\lambda_j)&=&\Biggl(\frac{1}{T}\sum
_{t=1}^Tf_tf_t'\Biggr)^{-1}\Biggl(\frac{1}{\sqrt T}\sum_{t=1}^T
f_te_{jt}\Biggr)+o_p(1)\\
&\xrightarrow{d}&\mathcal{N}(0,(\overline M_{f f})^{-1}\sigma
_j^2)
\end{eqnarray*}
and if $\sqrt N/T \to0$, then
\begin{eqnarray*}
\sqrt N(\hat f_t^{\mathrm{pc}}-f_t)&=&\Biggl(\frac{1}{N}\sum_{i=1}^N\lambda
_i\lambda
_i'\Biggr)^{-1}\Biggl(\frac1 {\sqrt N} \sum_{i=1}^N\lambda_ie_{it}
\Biggr)+o_p(1)\\
&\xrightarrow{d}&\mathcal{N}(0,(\overline M_{\lambda\lambda
})^{-1}\Upsilon(\overline M_{\lambda\lambda})^{-1}),
\end{eqnarray*}
where $\overline M_{\lambda\lambda}=\lim_{N\to\infty
}\frac1 N
\sum_{i=1}^N\lambda_i\lambda_i'$ and $\Upsilon=\lim_{N\to
\infty
}\frac1 N\sum_{i=1}^N\lambda_i\lambda_i'\sigma_i^2$.

The PC estimator in~\cite{r6} uses the identification restriction IC3.
Under IC3, we already show that the MLE satisfies, as $N, T\to\infty$,
\[
\sqrt T (\hat\lambda_j-\lambda_j)\xrightarrow{d}\mathcal{N}
(0,(\overline M_{f f})^{-1}\sigma_j^2).
\]
Theorem~\ref{thm7} above shows that, under IC3,
$\sqrt N (\hat f_t-f_t)\xrightarrow{d}\mathcal{N}(0, Q^{-1})$.
This result requires $\sqrt{N}/T\rightarrow0$, which is satisfied if
$N/T\rightarrow\Delta$.

While $\hat\lambda_j^{\mathrm{pc}}$ and $\hat\lambda_j$ have the same
limiting distribution, $\hat f_t^{\mathrm{pc}}$ is less efficient
than~$\hat f_t$. This follows because the sandwich form of the
covariance matrix $(\overline M_{\lambda\lambda})^{-1}\Upsilon
(\overline M_{\lambda\lambda })^{-1}$ is no smaller than $Q^{-1}$,
where $Q$ is the\vspace*{-1pt} limit of $\frac 1 N\sumiN\frac1 {\sigma_i^{2}}
\times\lambda_i \lambda_i'$. Moreover, under IC3, the MLE $\hat\lambda_j$
only requires $N, T\to\infty$, but the PC estimator
$\hat\lambda_j^{\mathrm{pc}}$ requires an additional assumption that
$\sqrt T/N \to0$. Furthermore, the maximum likelihood estimator
$\hat\lambda_j$ is consistent under fixed $N$, but $\hat
\lambda_j^{\mathrm{pc}}$ requires both $N$ and $T$ to be large in order
to have consistency. Of course, under fixed $N$, the limiting
distribution of MLE will have a~different (more complicated) asymptotic
covariance matrix; see~\cite{r17}.

To estimate\vspace*{-1pt} $\sigma_j^2$, the PC method would need to estimate the
individual residuals $\hat e_{it}=X_{it}-\hat\alpha_i -(\hat\lambda
_i^{\mathrm{pc}})' \hat f_t^{\mathrm{pc}}$ and then construct $\hat\sigma_j^2=\frac1 T
\sum_{t=1}^T\hat e_{jt}^2$. In case that $N$ is fixed, $f_t$ cannot be
consistently estimated, so $\hat e_{it}$ is inconsistent for $e_{it}$.
This further implies that $\hat\sigma_j^2$ is inconsistent for
$\sigma
_j^2$. In comparison, the MLE does not estimate the individuals $\hat
e_{it}$. The variances are estimated jointly with the factor loadings
$\lambda_j$ and with the matrix $M_{f f}$. The variance estimator
remains consistent under fixed~$N$.

Finally, the PC estimator for $\lambda_i$ satisfies (see~\cite{r6})
$\frac{1}{N}\sum_{i=1}^N\|\hat\lambda_i^{\mathrm{pc}}-\lambda_i\|
^2=O_p(\frac1
N)+O_p(\frac1 T)$,
while the MLE for $\lambda_i$ satisfies
$ \frac{1}{N}\sum_{i=1}^N\|\hat\lambda_i-\lambda_i \|^2=O_p(\frac
1 T)$.

\section{Computational issues}\label{sec8}
The maximum likelihood estimation can be
implemented via the EM algorithm and
is considered by~\cite{rubin1982algorithms}.
The EM algorithm is an iterated approach. To be specific, consider the
identification condition IC3. Once the estimator under IC3 is obtained,
estimators under other identification restrictions can be easily
obtained (to be discussed below). Under IC3, we only need to estimate
$\Lambda$ and $\Sigma_{ee}$ since $M_{f f}=I_r$.

Let $\theta^{(k)}=(\Lambda^{(k)},\Sigma_{ee}^{(k)})$ denote the
estimator at the $k$th iteration. The EM algorithm updates the
estimator according to
\begin{eqnarray*}
\Lambda^{(k+1)}&=&\Biggl[\frac1 T\sum_{t=1}^TE\bigl(z_tf_t'|Z, \theta
^{(k)}\bigr)\Biggr]\Biggl[\frac1 T\sum_{t=1}^TE\bigl(f_tf_t'|Z, \theta
^{(k)}\bigr)\Biggr]^{-1},
\\
\Sigma_{ee}^{(k+1)}&=& \diag\bigl(M_{z z}- \Lambda^{(k+1)}\Lambda
^{(k)\prime} \bigl(\Sigma_{z z}^{(k)}\bigr)^{-1} M_{z z}\bigr),
\end{eqnarray*}
where $\Sigma_{z z}^{(k)}=\Lambda^{(k)}{\Lambda^{(k)}}'+\Sigma
_{ee}^{(k)}$, and
\begin{eqnarray*}
\frac1 T\sum_{t=1}^TE\bigl(f_tf_t'|Z,\theta^{(k)}\bigr)&=&
{\Lambda^{(k)}}'\bigl(\Sigma_{z z}^{(k)}\bigr)^{-1}M_{z z}
\bigl(\Sigma_{z z}^{(k)}\bigr)^{-1}\Lambda^{(k)}\\
&&{}+I_r-{\Lambda^{(k)}}'\bigl(\Sigma_{z z}^{(k)}\bigr)^{-1}\Lambda^{(k)},
\\
\frac1 T\sum_{t=1}^TE\bigl(z_t f_t'|Z, \theta^{(k)}\bigr)&=&M_{z z} \bigl(\Sigma
_{z
z}^{(k)}\bigr)^{-1} \Lambda^{(k)}.
\end{eqnarray*}
This gives $\theta^{(k+1)}=(\Lambda^{(k+1)},\Sigma_{ee}^{(k+1)})$. The
iteration continues until $\|\theta^{(k+1)}-\theta^{(k)}\|$ is smaller
than a~preset tolerance. In the simulation reported below, we use the
principal components estimator as the starting value. Let $(\Lambda
^{\dag},\Sigma^\dag)$ denote the final round of iteration. Let
$\mathcal
{V}$ be the orthogonal matrix consisting of the eigenvectors of $\frac
1 N{\Lambda^\dag}'(\Sigma_{ee}^\dag)^{-1}\Lambda^\dag$
corresponding to
descending eigenvalues. Let $\hat\Lambda=\Lambda^\dag\mathcal{V}$ and
$\hat\Sigma_{ee}=\Sigma_{ee}^\dag$. Then $\hat\theta=(\hat
\Lambda,\hat
\Sigma_{ee})$ satisfies IC3. For general models,
\cite{wu1983convergence} shows that the EM solutions are stationary points
of the likelihood functions. For completeness, we provide a~direct and
simple proof of this claim for factor models in the supplement~\cite{bai2011Supplement}
(Section E).

It is interesting to note that, under large $N$ and large $T$, the
number of iterations needed to achieve convergence is smaller than
under either a~small~$N$ or a~small~$T$.
In Section E of the supplement~\cite{bai2011Supplement}, we also explain how to write a~computer
program so it runs fast.

Let $(\hat\Lambda, \hat\Sigma_{ee})$ denote the MLE under IC3. We
discuss how to obtain estimators that satisfy other identification restrictions.
First, note that $\hat\Sigma_{ee}$ is identical under IC1--IC5. We
only need to discuss how to obtain $\Lambda$ and~$M_{f f}$. Let $\hat
\Lambda^{\ell}$ and $\hat M_{f f}^{\ell}$ denote the MLE under
IC$\ell
$ $(\ell=1,\ldots,5)$. Let $\hat\Lambda_1$ denote the first $r\times r$
block of $\hat\Lambda$.
For IC1, let $\hat\Lambda^1=\hat\Lambda(\hat\Lambda_1)^{-1}$ and
$\hat
M_{f f}^1=\hat\Lambda_1\hat\Lambda_1'$. This new estimator
satisfies IC1.
For IC2, let $\hat\Lambda^2=\hat\Lambda(\frac1 N\hat\Lambda'\hat
\Sigma
_{ee}^{-1}\hat\Lambda)^{-1/2}$ and $\hat M_{f f}^2= \frac1 N\hat
\Lambda'\hat\Sigma_{ee}^{-1}\hat\Lambda$. Then this estimator
satisfies IC2.
For IC4, $M_{f f}=I_r$ is known. Let $\Lambda_1'=\mathcal{Q}\mathcal
{R}$ be the QR decomposition of $\Lambda_1'$ with $\mathcal{Q}$ an
orthogonal matrix and $\mathcal{R}$ an upper triangular matrix. Define
$\hat\Lambda^4=\hat\Lambda\mathcal{Q}$. Then $\hat\Lambda^4$
satisfies IC4.
Finally, consider IC5. Let $\mathcal{W}$ be the diagonal matrix with
its diagonal elements the same as the first $r\times r$ block of $\hat
\Lambda^4$.
Let $\hat\Lambda^5=\hat\Lambda^4\mathcal{W}^{-1}$ and $ \hat M_{f
f}^5=\mathcal{W}{\mathcal{W}}'$. Then IC5 is satisfied.\vspace*{1pt}

%
\begin{table}
\caption{The performance of MLE and PC}\label{table2}
\begin{tabular*}{\tablewidth}{@{\extracolsep{\fill}}lccccccc@{}}
\hline
& & \multicolumn{3}{c}{\textbf{MLE}} & \multicolumn{3}{c@{}}{\textbf{PC}}\\[-4pt]
& & \multicolumn{3}{c}{\hrulefill} & \multicolumn{3}{c@{}}{\hrulefill}\\
$\bolds{N}$ & $\bolds{T}$ & $\bolds{\Lambda}$ & $\bolds{F}$
& $\bolds{\Sigma_{ee}}$ & $\bolds{\Lambda}$ & $\bolds{F}$ &
$\bolds{\Sigma_{ee}}$\\
\hline
\hphantom{0}$10$ & \hphantom{0}$30$ &0.4818&0.3473 & 0.8432& 0.4058&0.2744 &0.7991 \\
\hphantom{0}$30$ & \hphantom{0}$30$ &0.7276 &0.7995 &0.9273 &0.6391 &0.6450 & 0.9223 \\
\hphantom{0}$50$ & \hphantom{0}$30$ & 0.7676& 0.8973& 0.9303& 0.7221& 0.7953 & 0.9302 \\
$100$ & \hphantom{0}$30$ &0.7874 &0.9555 &0.9308 & 0.7679& 0.9006& 0.9312 \\
$150$ & \hphantom{0}$30$ &0.7941 & 0.9719&0.9310 & 0.7823&0.9347 & 0.9315 \\
[4pt]
\hphantom{0}$10$ & \hphantom{0}$50$ &0.6080 &0.4153 & 0.8951&0.4875&0.2975 &0.8187 \\
\hphantom{0}$30$ & \hphantom{0}$50$ &0.8383 &0.8407 &0.9583 &0.7751 &0.7113 & 0.9499 \\
\hphantom{0}$50$ & \hphantom{0}$50$ &0.8589 &0.9161 &0.9590 & 0.8306&0.8341 & 0.9569 \\
$100$ & \hphantom{0}$50$ &0.8722 &0.9624 &0.9592 &0.8613 &0.9198 & 0.9591 \\
$150$ & \hphantom{0}$50$ & 0.8764&0.9764 & 0.9592&0.8697 &0.9475 & 0.9593 \\
[4pt]
\hphantom{0}$10$ & $100$ &0.7563&0.4939 & 0.9448& 0.5878 &0.3298 &0.8345 \\
\hphantom{0}$30$ & $100$ &0.9182 & 0.8614&0.9793 & 0.8789&0.7519 &0.9700 \\
\hphantom{0}$50$ & $100$ &0.9292 &0.9245 &0.9798 &0.9135 &0.8572 &0.9770 \\
$100$ & $100$ &0.9362 &0.9668 &0.9798 & 0.9305&0.9308 & 0.9792 \\
$150$ & $100$ &0.9383& 0.9788 &0.9799 & 0.9349& 0.9545 & 0.9798 \\
\hline
\end{tabular*}
\end{table}

We now consider the finite sample properties of the MLE.
Data are generated according to $z_{it}=\lambda_i' f_t+e_{it}$ with
$r=2$, where $\lambda_i, f_t$ are i.i.d. $\mathcal{N}(0, I_2)$ and
$e_{it}$ follows $\mathcal{N}(0,\sigma_i^2)$ with $\sigma
_i^2=0.1+10\times U_i$, and $U_i$ are i.i.d. uniform on $[0,1]$. Adding
0.1 to the variance avoids near-zero values. We consider combinations
of $T=30,50,100$ and $N=10,30,50,100,150$. Estimators under different
identification conditions only differ up to a~rotation matrix, so only
IC3 will be considered. We also compute the principal components (PC)
estimator for comparison.
To measure the accuracy between $\hat\Lambda$ and $\Lambda$ (both are
$N\times2)$, we compute the second (the smallest nonzero) canonical
correlation between them. Canonical correlation is widely used as a~measure of goodness-of-fit in factor analysis; see, for example,
\cite{r13} and~\cite{goyal-et-al}. Similarly, we also compute the second
canonical correlation between $\hat F$ and $F$. For the estimated
variances, we calculate the squared correlation between $\operatorname
{diag}(\hat\Sigma_{ee})$ and $\operatorname{diag}(\Sigma_{ee})$. The
corresponding values for the principal components estimators are also
computed. Table~\ref{table2} reports the average canonical correlations based on
5000 repetitions for each $(N,T)$ combination.

The results\vspace*{1pt} suggest that the precision of $\hat\Lambda$ is closely
tied to the size of $T$ and the precision of $\hat F$ is tied to $N$.
This is consistent with the theory. For all $(N,T)$ combinations, the
MLE dominates PC. The domination becomes less important for $N\ge50$
and $T\ge50$ for estimating factor loadings. But for small $N$, no
matter how large is $T$, MLE noticeably outperforms PC.
For the estimated factors, there is still noticeable outperformance
even under large $N$ and $T$. These are all consistent with the theory.

\section{Conclusion}\label{sec9}
In this paper we have developed an inferential theory for factor models
of high dimension.
We study the maximum likelihood estimator under five different sets of
identification restrictions.
Consistency, rate of convergence and the limiting distributions are
derived. Unlike the principal component methods, the estimators are
shown to be efficient under the model assumptions.
While both the factor loadings and factors are treated as parameters
(nonrandom), the key to efficiency is not to simultaneously estimate
both the factor loadings and the factors. If
$N$ is relatively small compared with $T$, the efficient approach is
to estimate the individual factor loadings $(\lambda_i)$ and the sample
moment of the factor scores ($f_t$), not the individual scores. The
sample moment contains only $r(r+1)/2$ unknown elements. If the factor
scores $f_t$ are of interest, they can be estimated by the generalized
least squares in a~separate stage. The estimated factor scores are also
shown to be efficient under the model assumptions. The opposite
procedure should be adopted\vadjust{\goodbreak} if $N$ is much larger than $T$. In the
latter case, we estimate the individual factor scores and the sample
moment of the factor loadings. If $N$ and $T$ are comparable, the
choice of procedures boils down which heteroscedasticity,
cross-sectional dimension or the time dimension, is the object of
interest. The paper also provides a~novel approach to consistency in
the presence of a~large and increasing number of parameters.\vspace*{-3pt}

%
%
\begin{appendix}\label{app}
\section*{\texorpdfstring{Appendix: Proof of Proposition \lowercase{\protect\ref{thm1}}}{Appendix: Proof of Proposition 5.1}}
\vspace*{-3pt}

The following notation will be used throughout:
\begin{eqnarray*}
\hat H &=&(\hat\Lambda'\hat\Sigma_{ee}^{-1}\hat\Lambda
)^{-1}, \\[-2pt]
\hat H_N &=& N\cdot\hat H=( N^{-1}\hat\Lambda'\hat\Sigma
_{ee}^{-1}\hat
\Lambda)^{-1}, \\[-2pt]
\hat G&=&(\hat M_{f f}^{-1}+\hat\Lambda'\hat\Sigma_{ee}^{-1}\hat
\Lambda
)^{-1},\\[-2pt]
\hat G_N &=& N \cdot\hat G, \\[-2pt]
\xi_t&=&(e_{1t},e_{2t},\ldots, e_{rt})'.
\end{eqnarray*}
From $(A+B)^{-1}=A^{-1}-A^{-1}B(A+B)^{-1}$, we have
$\hat H=\hat G(I - \hat M_{f f}^{-1}\hat G)^{-1}$. From $\Sigma_{z
z}=\Lambda M_{f f}\Lambda'+\Sigma_{ee}$, we have
%
%
\begin{equation}\label{inv}
\Sigma_{z z}^{-1}=\Sigma_{ee}^{-1}-\Sigma_{ee}^{-1}\Lambda(M_{f
f}^{-1}+\Lambda'\Sigma_{ee}^{-1}\Lambda)^{-1}\Lambda'\Sigma_{ee}^{-1}.
\end{equation}
It follows %
%
%
\begin{eqnarray}\label{seca2}
\hat\Lambda'\hat\Sigma_{z z}^{-1}&=&\hat\Lambda'\hat\Sigma
_{ee}^{-1}-\hat\Lambda'\hat\Sigma_{ee}^{-1}\hat\Lambda(\hat M_{f
f}^{-1}+\hat\Lambda'\hat\Sigma_{ee}^{-1}\hat\Lambda)^{-1}\hat
\Lambda
'\hat\Sigma_{ee}^{-1}\nonumber\\[-9pt]\\[-9pt]
&=&\hat M_{f f}^{-1}\hat G\hat\Lambda'\hat
\Sigma_{ee}^{-1}.\nonumber
\end{eqnarray}

To prove Proposition~\ref{thm1}, we use a~superscript ``$^*$'' to denote
the true parameters, for example, $\Lambda^*$, $\Sigma_{ee}^*$,
$f_t^*$, etc. The variables without the superscript ``$^*$'' denote the
function arguments (input variables) in the likelihood function.

Let $\theta=(\lambda_1,\ldots,\lambda_n,\sigma_1^2,\ldots,\sigma_n^2,
M_{f f})$ and
let $\Theta$ be a~parameter set such that $C^{-2}\le\sigma_t^2\le C^2$,
$M_{f f}$ is positive definite matrices with elements bounded. We assume
$\theta^*=(\lambda_1^*,\ldots,\lambda_n^*,\sigma_1^{2*},\ldots
,\sigma_n^{2*}$
$, M_{f f}^*)$
is an interior point of $\Theta$. For simplicity, we also write
$\theta
=(\Lambda, \Sigma_{ee},
M_{f f})$ and $\theta^*=(\Lambda^*, \Sigma_{ee}^*, M_{f f}^*)$.
\begin{pf*}{Proof of Proposition~\ref{thm1}}
The centered likelihood function can be written as
\[
L(\theta)=\overline L(\theta)+R(\theta),
\]
where
\[
\overline L(\theta) = - {\frac{1}{N}\ln}|\Sigma_{z z} | - \frac
{1}{N}\tr( \Sigma_{z z}(\theta^*)\Sigma_{z z}^{ - 1} ) +
1 +
{\frac{1}{N}\ln}|\Sigma(\theta^*)|
\]
and $R(\theta)=-\frac{1}{N}\tr((M_{z z}-\Sigma_{z z}(\theta
^*))\Sigma_{z z}^{-1})$.
Note that $1 + \frac{1}{N}{\sum_{i = 1}^N \ln}|\Sigma(\theta^*)|$ does
not depend on any unknown parameters and is for the purpose of centering.

Lemma A.2~\cite{bai2011Supplement} implies that $\sup_{\theta}|R(\theta)|=o_p(1)$. In
particular, we have
$|R(\hat\theta)|=o_p(1)$ and\vadjust{\goodbreak} $|R(\theta^*)|=o_p(1)$. So $|R(\theta
^*)-R(\hat\theta)|=o_p(1)$. Since $\hat\theta$ maximizes $L(\theta)$,
it follows $\overline L(\hat\theta)+R(\hat\theta))\ge\overline
L(\theta
^*)+R(\theta^*)$. Hence we have $\overline L(\hat\theta)\ge
\overline
L(\theta^*)+R(\theta^*)-R(\hat\theta)\ge\overline L(\theta
^*)-|o_p(1)|$. However, the function $\overline L(\theta)$ achieves its
maximum at $\theta^*$, so $\overline L(\hat\theta)\le\overline
L(\theta
^*)$. Since $\overline L(\theta^*)$ is normalized to zero, we have
$\overline L(\hat\theta)\ge-|o_p(1)|$ and $\overline L(\hat\theta
)\le0$. It follows that $\overline L(\hat\theta)=o_p(1)$.

Notice $|\Sigma_{z z}|=|\Sigma_{ee}|\cdot|I_r+M_{f f}\Lambda
'\Sigma
_{ee}^{-1}\Lambda|$. But
$|I_r+M_{f f}\Lambda'\Sigma_{ee}^{-1}\Lambda|=O(N)$. Similarly
$|\Sigma_{z z}(\theta^*)|=|\Sigma_{ee}^*|\cdot|I_r+M_{f
f}^*\Lambda^{*\prime}
\Sigma_{ee}^{*-1}\Lambda^*|$, thus uniformly on $\Theta$,
\[
- {\frac{1}{N}\ln}|\Sigma_{z z}|
+{\frac1 N \ln}|\Sigma_{z z}(\theta^*)|= - {\frac{1}{N}\ln}|\Sigma_{ee}|
+{\frac1 N \ln}|\Sigma_{ee}^*| + O\biggl(\frac{\ln(N)} N \biggr).
\]

Next, from $\Sigma_{z z}(\theta^*)=\Lambda^{*} M_{f f}^* \Lambda
^{*\prime}+\Sigma_{ee}^*$, we have $\Sigma_{z z}(\theta^*)\Sigma_{z
z}^{-1}=\break\Lambda^{*} M_{f f}^*
\Lambda^{*\prime}\times\Sigma_{z z}^{-1} +\Sigma_{ee}^*
\Sigma_{z z}^{-1}$. Using the formula for $\Sigma_{z z}^{-1}$, we have
$ \tr(\Sigma_{ee}^* \Sigma_{z z}^{-1})=\tr(\Sigma_{ee}^*\Sigma
_{ee}^{-1}) +O(1)$, because $\tr[\Sigma_{ee}^*\Sigma _{ee}^{-1}\Lambda
(M_{f f}^{-1}+\Lambda'\Sigma_{ee}^{- 1} \Lambda)^{-1}\Lambda '\Sigma
_{ee}^{-1}]=O(1)$. The latter follows since the matrix in the square
bracket is bounded in norm by
$C^4\|\Lambda'\Sigma_{ee}^{-1}\Lambda\times
(M_{f f}^{-1}+\Lambda'\Sigma_{ee}^{- 1}\Lambda)^{-1}\| \leq C^4
\|I_r\|$ due\vspace*{1pt} to the bound on $\sigma_i^2$ and~$\sigma_i^{*2}$. Thus divided by $N$, we have
\[
\frac1 N \tr[\Sigma_{z z}(\theta^*)\Sigma_{z z}^{-1}]= \frac1 N
\tr
[ \Lambda^{*} M_{f f}^* \Lambda^{*\prime}\Sigma_{z z}^{-1}] + \frac1 N
\tr
(\Sigma_{ee}^*\Sigma_{ee}^{-1}) + O\biggl(\frac1 N\biggr).
\]
Notice ${\ln}|\Sigma_{ee}|=\sumiN\ln\sigma_i^2$ and $\tr(\Sigma
_{ee}^*\Sigma_{ee}^{-1})=\sumiN\sigma_i^{*2}/\sigma_i^2$; we have
proved that
\[
\overline L(\theta)=-\frac{1}{N}\sum_{i=1}^N \biggl(\ln\sigma_i^2 +
\frac
{\sigma_i^{*2} }{\sigma_i^2 }- 1 - \ln\sigma_i^{*2}\biggr) - \frac
{1}{N}\tr(\Lambda^*M_{f f}^*\Lambda^{*\prime}\Sigma_{z z}^{-1}) +O
\biggl(\frac{\ln N} N\biggr)
\]
uniformly on $\Theta$. By $\overline L(\hat\theta)=o_p(1)$, it
follows that
\[
-\frac{1}{N}\sum_{i = 1}^N \biggl(\ln\hat\sigma_i^2 + \frac{\sigma
_i^{*2}}{\hat\sigma_i^2} - 1 - \ln\sigma_i^{*2}\biggr) - \frac{1}{N}\tr
(\Lambda^* M_{f f}^* \Lambda^{*\prime} \hat\Sigma_{z
z}^{-1})\xrightarrow{p}0.
\]
A~key observation is that both terms are nonpositive; it follows
%
%
\begin{eqnarray}\label{pthm4}
\frac{1}{N}\sum_{i = 1}^N\biggl(\ln\hat\sigma_i^2+\frac{\sigma
_i^{*2}}{\hat
\sigma_i^2} - 1 -\ln\sigma_i^{*2} \biggr)&\xrightarrow{p}& 0,
\\
%
%
\label{pthm5}
\frac{1}{N}\tr(\Lambda^*M_{f f}^*\Lambda^{*\prime}\hat\Sigma_{z
z}^{-1})&\xrightarrow{p}&0.
\end{eqnarray}
Consider the function $f(x)=\ln x+\frac{\sigma_i^{*2}}{x}-\ln\sigma
_i^{*2}-1$. Given that $0< C^{-2}\le\sigma_i^2\le C^2 < \infty$ for
$C>1$, for any $x\in[C^{-2},C^2]$, there exists a~constant $b$ (e.g.,
take $b=\frac{1}{4C^4}$), such that $f(x)\ge b(x-\sigma_i^{*2})^2$.
It follows
\[
o_p(1) = \frac{1}{N}\sum_{i = 1}^N\biggl(\ln\hat\sigma_i^2 +\frac
{\sigma
_i^{*2}}{\hat\sigma_i^2} - 1 - \ln\sigma_i^{*2}\biggr) \ge b\frac
{1}{N}\sum
_{i = 1}^N(\hat\sigma_i^2-\sigma_i^{*2})^2.
\]
The above implies
%
%
\begin{equation}\label{convar}
\frac{1}{N}\sum_{i = 1}^N (\hat\sigma_i^2- \sigma
_i^{*2})^2\xrightarrow{p}0.
\end{equation}
Now we turn to~(\ref{pthm5}). By~(\ref{inv}), we have
\begin{eqnarray*}
&&
\frac{1}{N}\tr(\Lambda^* M_{f f}^* \Lambda^{*\prime} \hat\Sigma_{z
z}^{ -
1})\\
&&\qquad = \frac{1}{N}\tr\bigl(M_{f f}^* \Lambda^{*\prime} [\hat\Sigma
_{ee}^{-1}-\hat\Sigma_{ee}^{-1}\hat\Lambda(\hat M_{f f}^{-1}+\hat
\Lambda'\hat\Sigma_{ee}^{-1}\hat\Lambda)^{- 1} \hat\Lambda' \hat
\Sigma
_{ee}^{-1}]\Lambda^*\bigr).
\end{eqnarray*}
From $ (\hat M_{f f}^{-1}+\hat\Lambda'\hat\Sigma_{ee}^{-1}\hat
\Lambda
)^{-1} = (\hat\Lambda'\hat\Sigma_{ee}^{ - 1} \hat\Lambda)^{-1} -
(\hat\Lambda'\hat\Sigma_{ee}^{ - 1} \hat\Lambda)^{-1} \hat M_{f
f}^{-1}(\hat M_{f f}^{ - 1} +\break \hat\Lambda'\times
\hat\Sigma_{ee}^{ - 1}
\hat\Lambda)^{-1}$, we obtain
\begin{eqnarray*}
&&\frac{1}{N}\tr(\Lambda^* M_{f f}^* \Lambda^{*\prime} \hat\Sigma_{z z}^{
- 1})\\
&&\qquad= \frac{1}{N}\tr[
M_{f f}^* \Lambda^{*\prime} \hat\Sigma_{ee}^{ - 1} \Lambda^* - M_{f
f}^* \Lambda^{*\prime} \hat\Sigma_{ee}^{ - 1} \hat\Lambda(\hat\Lambda'
\hat\Sigma_{ee}^{ - 1} \hat\Lambda)^{ - 1} \hat\Lambda' \hat
\Sigma
_{ee}^{ - 1} \Lambda^* ]
\\
&&\qquad\quad{} +
\tr[ M_{f f}^* \Lambda^{*\prime} \hat\Sigma_{ee}^{ - 1} \hat
\Lambda
( \hat\Lambda' \hat\Sigma_{ee}^{ - 1} \hat\Lambda)^{ - 1} \hat
M_{f f}^{ - 1} ( \hat M_{f f}^{ - 1} + \hat\Lambda' \hat\Sigma
_{ee}^{ - 1} \hat\Lambda)^{ - 1} \hat\Lambda' \hat\Sigma_{ee}^{ -
1} \Lambda^* ].
\end{eqnarray*}
Both expressions are nonnegative; by~(\ref{pthm5}), we must have
%
%
\begin{equation}\label{pthm7}
\frac{1}{N}\tr\bigl(M_{f f}^* \Lambda^{*\prime}\hat\Sigma_{ee}^{ - 1}
\Lambda
^* - M_{f f}^* \Lambda^{*\prime}\hat\Sigma_{ee}^{ - 1} \hat\Lambda( \hat
\Lambda'\hat\Sigma_{ee}^{ - 1} \hat\Lambda)^{ - 1} \hat\Lambda
'\hat\Sigma_{ee}^{ - 1} \Lambda^* \bigr)\xrightarrow{p}0\hspace*{-12pt}
\end{equation}
and
%
%
\begin{eqnarray}\label{pthm8}
&&\frac{1}{N}\tr\bigl( M_{f f}^* \Lambda^{*\prime}\hat\Sigma_{ee}^{ - 1}
\hat
\Lambda( \hat\Lambda'\hat\Sigma_{ee}^{ - 1} \hat\Lambda)^{ - 1}
\nonumber\\[-8pt]\\[-8pt]
&&\qquad{}\times\hat M_{f f}^{ - 1}( \hat M_{f f}^{ - 1} + \hat\Lambda'\hat\Sigma
_{ee}^{ - 1} \hat\Lambda)^{ - 1} \hat\Lambda'\hat\Sigma_{ee}^{ -
1} \Lambda^* \bigr)\xrightarrow{p}0.\nonumber
\end{eqnarray}
By~(\ref{convar}) and Lemma A.4~\cite{bai2011Supplement}, $\frac{1}{N}\tr( M_{f f}^* \Lambda
^{*\prime}\hat\Sigma_{ee}^{ - 1} \Lambda^*)\xrightarrow{p}C^*> 0$, say. From
(\ref{pthm7})
\[
\frac{1}{N}\tr( {M_{f f}^* \Lambda^{*\prime} \hat\Sigma_{ee}^{ - 1}
\hat
\Lambda( {\hat\Lambda' \hat\Sigma_{ee}^{ - 1} \hat\Lambda})^{ - 1}
\hat\Lambda'\hat\Sigma_{ee}^{ - 1} \Lambda^* } ) \xrightarrow
{p}C^*> 0
\]
with the same $C^*$.
The preceding result and~(\ref{pthm8}) imply
\[
\hat M_{f f}^{ - 1} ( {\hat M_{f f}^{ - 1} + \hat\Lambda'\hat
\Sigma
_{ee}^{ - 1} \hat\Lambda})^{ - 1}=o_p(1).
\]
By assumption, we confine $M_{f f}$ on a~compact set, that is, $\hat
M_{f f}=O_p(1)$. By the definition of $\hat G$, we have $\hat
G=o_p(1)$. From $\hat H=\hat G(I-\hat M_{f f}^{-1}\hat G)^{-1}$, we
have $\hat H=o_p(1)$. We obtain the following result:
%
%
\begin{equation}\label{pthm9}
\hat G=o_p(1);\qquad \hat H=o_p(1).
\end{equation}
The matrix on the left-hand side of~(\ref{pthm7}) is semi-positive
definite and is finite dimensional ($r\times r)$, its trace is $o_p(1)$
if and only if every entry is~$o_p(1)$.
Thus we have
%
\[
\frac{1}{N}\bigl(M_{f f}^* \Lambda^{*\prime} \hat\Sigma_{ee}^{ - 1}
\Lambda^* - M_{f f}^* \Lambda^{*\prime} \hat\Sigma_{ee}^{ - 1} \hat
\Lambda( {\hat\Lambda'\hat\Sigma_{ee}^{ - 1} \hat\Lambda})^{ - 1}
\hat\Lambda'\hat\Sigma_{ee}^{ - 1} \Lambda^*\bigr)\xrightarrow{p}0.
\]
Pre-multiplying both sides by $M_{f f}^{*-1}$ gives
%
%
\begin{equation}\label{pthm10}
\frac{1}{N}\Lambda^{*\prime} \hat\Sigma_{ee}^{ - 1} \Lambda^* - \frac1 N
\Lambda^{*\prime} \hat\Sigma_{ee}^{ - 1} \hat\Lambda(\hat\Lambda'\hat
\Sigma_{ee}^{ - 1} \hat\Lambda)^{ - 1} \hat\Lambda'\hat\Sigma
_{ee}^{ - 1} \Lambda^*\xrightarrow{p} 0.
\end{equation}
The second term on the left-hand side can be rewritten as
\[
[\Lambda^{*\prime} \hat\Sigma_{ee}^{ - 1}
\hat\Lambda( {\hat\Lambda'\hat\Sigma_{ee}^{ - 1} \hat\Lambda})^{
- 1} ]\biggl( \frac1 N {\hat\Lambda'\hat\Sigma_{ee}^{ - 1}
\hat
\Lambda} \biggr)[( \hat\Lambda'\hat\Sigma_{ee}^{ - 1} \hat
\Lambda)^{-1}\hat\Lambda'\hat\Sigma_{ee}^{ - 1} \Lambda^*].
\]
Let $A\equiv( {\hat\Lambda- \Lambda^* })' \hat\Sigma_{ee}^{ - 1}
\hat\Lambda\hat H$, where $\hat H=(\hat\Lambda'\hat\Sigma
_{ee}^{-1}\hat\Lambda)^{-1}$. It follows that $\Lambda^{*\prime} \hat
\Sigma
_{ee}^{ - 1} \times\hat\Lambda (\hat\Lambda'\hat\Sigma_{ee}^{ - 1}
\hat\Lambda)^{ - 1}=(I_r-A)$. So~(\ref{pthm10}) is equivalent to
\[
\frac{1}{N}\Lambda^{*\prime}\hat\Sigma_{ee}^{ - 1} \Lambda^*-(I_r -
A)\frac
{1}{N}\hat\Lambda'\hat\Sigma_{ee}^{ - 1} \hat\Lambda(I_r - A)'
\xrightarrow{p}0.
\]
However, $\frac{1}{N}\Lambda^{*\prime}\hat\Sigma_{ee}^{ - 1} \Lambda^*
=\frac{1}{N}\Lambda^{*\prime}\Sigma_{ee}^{*- 1} \Lambda^*+o_p(1)$ by Lemma
A.4~\cite{bai2011Supplement} and~(\ref{convar}), thus
%
%
\begin{equation}\label{pthm17}
\frac{1}{N}\Lambda^{*\prime} \Sigma_{ee}^{ *- 1} \Lambda^*-(I_r -
A)\frac
{1}{N}\hat\Lambda'\hat\Sigma_{ee}^{ - 1} \hat\Lambda(I_r - A)'
\xrightarrow{p}0.
\end{equation}
Because the first term is of full rank in the limit, the second term is
also of full rank. This implies that $I_r-A$ in the limit is of full rank.

Meanwhile, equation~(\ref{pthm10}) can be expressed alternatively as
\begin{eqnarray*}
&&\frac{1}{N}(\hat\Lambda- \Lambda^*)'\hat\Sigma_{ee}^{ - 1}(\hat
\Lambda- \Lambda^*)\\
&&\qquad{} - \frac{1}{N}(\hat\Lambda- \Lambda^*)' \hat
\Sigma_{ee}^{ - 1} \hat\Lambda\biggl(\frac{1}{N}\hat\Lambda' \hat
\Sigma_{ee}^{ - 1} \hat\Lambda\biggr)^{ - 1} \frac{1}{N}\hat\Lambda'
\hat\Sigma_{ee}^{ - 1} (\hat\Lambda- \Lambda^*)\xrightarrow{p}0,
\end{eqnarray*}
which can also be written as, in terms of $A$,
%
%
\begin{equation}\label{pthm12}
\frac{1}{N}( \hat\Lambda- \Lambda
^*)'\hat\Sigma_{ee}^{ - 1}( {\hat\Lambda- \Lambda^* })- A~\biggl(\frac
{1}{N}\hat\Lambda' \hat\Sigma_{ee}^{ - 1}\hat\Lambda
\biggr)A'\xrightarrow{p}0.
\end{equation}
Both~(\ref{pthm17}) and~(\ref{pthm12}) will be useful in establishing
consistency.

We now make use of the first-order conditions.
The first-order condition~(\ref{sec271}), by~(\ref{seca2}), can be
simplified as $\hat\Lambda'\hat\Sigma_{ee}^{-1}(M_{z z}-\hat
\Sigma_{z
z})=0$. This gives
\begin{eqnarray*}
&&\hat\Lambda'\hat\Sigma_{ee}^{-1}\Biggl(\Lambda^* M_{f f}^*\Lambda
^{*\prime}+\Lambda^*\frac{1}{T}\sum_{t=1}^Tf_t^*e_t'+\frac{1}{T}\sum
_{t=1}^Te_tf_t^{*\prime}\Lambda^{*\prime}
\\
&&\hspace*{17pt}\qquad{}
+\frac{1}{T}\sum_{t=1}^T(e_te_t'-\Sigma_{ee}^*)
+\Sigma_{ee}^*-\hat\Lambda\hat M_{f f}\hat\Lambda'-\hat\Sigma
_{ee}\Biggr)=0.
\end{eqnarray*}
For simplicity, we neglect the smaller-order term $\hat\Lambda'\hat
\Sigma_{ee}^{-1}\bar e\bar e'$. The $j$th column of the above equation
can be written as (after some algebra),
%
%
\begin{eqnarray}\label{pthm20}
\hat\lambda_j-\lambda_j^*&=& -\hat M_{f f}^{-1}(\hat M_{f f}- M_{f
f}^*)\lambda_j^*\nonumber\\
&&{}-\hat M_{f f}^{-1}\hat H\hat\Lambda'\hat\Sigma
_{ee}^{-1}(\hat\Lambda-\Lambda^*)M_{f f}^*\lambda_j^*
\nonumber\\
&&{}
+\hat M_{f f}^{-1}\hat H\hat\Lambda
'\hat\Sigma_{ee}^{ - 1} \Lambda^*\Biggl(\frac{1}{T}\sum_{t = 1}^T
f_t^*e_{jt}\Biggr)\\
&&{}+\hat M_{f f}^{-1} \hat H\hat\Lambda'\hat\Sigma
_{ee}^{-1}\Biggl(\frac{1}{T}\sum_{t=1}^Te_tf_t^{*\prime}\Biggr)\lambda_j^*
- \hat M_{f f}^{ - 1} \hat H\hat\lambda_j \frac{1}{\hat\sigma
_j^2}(\hat\sigma_j^2-\sigma_j^{*2})\nonumber\\
&&{}+ \hat M_{f f}^{ - 1} \hat H
\Biggl(\sum_{i = 1}^N \frac{1}{\hat\sigma_i^2}\hat\lambda_i\frac
{1}{T}\sum
_{t=1}^T [e_{it}e_{jt}- E(e_{it} e_{jt})] \Biggr).\nonumber
\end{eqnarray}
Consider the first-order condition~(\ref{sec273}). By the method
analogous to the one in deducing~(\ref{pthm20}), we have
%
%
\begin{eqnarray} \label{pthm21}
\hat M_{f f}- M_{f f}^*&=&-\hat H\hat\Lambda'\hat\Sigma
_{ee}^{-1}(\hat
\Lambda-\Lambda^*)M_{f f}^*-M_{f f}^* (\hat\Lambda- \Lambda
^*)'\hat
\Sigma_{ee}^{-1}\hat\Lambda\hat H
\nonumber\\
&&{}
+\hat H\hat\Lambda'\hat\Sigma_{ee}^{-1}(\hat\Lambda-\Lambda
^*)M_{f
f}^*(\hat\Lambda-\Lambda^*)'\hat\Sigma_{ee}^{- 1}\hat\Lambda\hat H
\nonumber\\
&&{}+\hat H\hat\Lambda' \hat\Sigma_{ee}^{-1}\Lambda^*\Biggl(\frac{1}{T}\sum
_{t= 1}^T f_t^*e_t'\Biggr)\hat\Sigma_{ee}^{-1} \hat\Lambda\hat H
\nonumber\\[-8pt]\\[-8pt]
&&{} + \hat
H\hat\Lambda'\hat\Sigma_{ee}^{ - 1}\Biggl(\frac{1}{T}\sum_{t = 1}^Te_t
f_t^{*\prime}\Biggr)\Lambda^{*\prime}\hat\Sigma_{ee}^{-1}\hat\Lambda\hat H
\nonumber\\
&&{}+ \hat H\Biggl(\sum_{i = 1}^N\sum_{j = 1}^N \frac{1}{\hat\sigma_i^2
\hat\sigma_j^2}\hat\lambda_i \hat\lambda_j' \frac{1}{T}\sum_{t =
1}^T [e_{it} e_{jt}-E(e_{it} e_{jt})]\Biggr)\hat H \nonumber\\
&&{}- \hat H\sum_{i =
1}^N\frac{1}{\hat\sigma_i^4 }\hat\lambda_i \hat\lambda_i'( \hat
\sigma_i^2 - \sigma_i^{*2}) \hat H.\nonumber
\end{eqnarray}
Substituting~(\ref{pthm21}) into~(\ref{pthm20}), we obtain
%
%
\begin{eqnarray}\label{pthm22}\quad
\hat\lambda_j-\lambda_j^*
&=& \hat M_{f f}^{-1} M_{f f}^* (\hat
\Lambda
-\Lambda^*)'\hat\Sigma_{ee}^{-1}\hat\Lambda\hat H\lambda_j^*\nonumber\\
&&{}-\hat
M_{f f}^{-1}\hat H\hat\Lambda'\hat\Sigma_{ee}^{-1}(\hat\Lambda
-\Lambda
^*)M_{f f}^* (\hat\Lambda- \Lambda^*)' \hat\Sigma_{ee}^{-1}\hat
\Lambda
\hat H\lambda_j^*
\nonumber\\
&&{}
-\hat M_{f f}^{-1}\hat H\hat\Lambda' \hat\Sigma_{ee}^{-1} \Lambda
^*\Biggl(\frac{1}{T}\sum_{t=1}^T f_t^*e_t'\Biggr)\hat\Sigma_{ee}^{-1}\hat
\Lambda
\hat H\lambda_j^* \nonumber\\
&&{}-\hat M_{f f}^{-1}\hat H\hat\Lambda'\hat\Sigma
_{ee}^{-1}\Biggl(\frac{1}{T}\sum_{t = 1}^T e_t f_t^{*\prime}\Biggr)\Lambda^{*\prime} \hat
\Sigma_{ee}^{-1}\hat\Lambda\hat H\lambda_j^*
\nonumber\\
&&{}- \hat M_{f f}^{ - 1} \hat H\Biggl( {\sum_{i = 1}^N {\sum_{j = 1}^N
{\frac{1}{{\hat\sigma_i^2 \hat\sigma_j^2 }}\hat\lambda_i \hat
\lambda_j' \frac{1}{T}} } \sum_{t = 1}^T [e_{it} e_{jt} - E( {e_{it}
e_{jt} })] } \Biggr)\hat H\lambda_j^*
\nonumber\\
&&{}
+ \hat M_{f f}^{ - 1} \hat H\sum_{i
= 1}^N \frac{1}{\hat\sigma_i^4}\hat\lambda_i \hat\lambda_i'
(\hat
\sigma_i^2 - \sigma_i^{*2} )\hat H\lambda_j ^*\nonumber\\
&&{}+ \hat M_{f
f}^{-1}\hat H\hat\Lambda'\hat\Sigma_{ee}^{-1}\Biggl(\frac{1}{T}\sum
_{t=1}^Te_t f_t^{*\prime}\Biggr)\lambda_j^*
\nonumber\\
&&{}+ \hat M_{f f}^{-1}\hat H\hat\Lambda'\hat\Sigma_{ee}^{-1}\Lambda
^*\Biggl(\frac{1}{T}\sum_{t = 1}^T f_t^*e_{jt}\Biggr)\\
&&{}+ \hat M_{f f}^{ - 1} \hat
H\Biggl( {\sum_{i = 1}^N {\frac{1}{{\hat\sigma_i^2 }}\hat\lambda_i
\frac{1}{T}\sum_{t = 1}^T [e_{it} e_{jt} - E( {e_{it} e_{jt} })] }
}\Biggr)
\nonumber\\
&&{}- \hat M_{f f}^{- 1}\hat H\hat\lambda_j \frac{1}{\hat\sigma_j^2 }(
\hat\sigma_j^2-\sigma_j^{*2}).\nonumber
\end{eqnarray}
Consider~(\ref{pthm21}). The fifth term of the right-hand side can be
written as
\[
\hat H\hat\Lambda'\hat\Sigma_{ee}^{ - 1}\Biggl(\frac{1}{T}\sum_{t=1}^Te_t
f_t^{*\prime}\Biggr)-\hat H\hat\Lambda'\hat\Sigma_{ee}^{-1}\Biggl(\frac{1}{T}\sum
_{t=1}^T e_tf_t^{*\prime}\Biggr)A,
\]
where $A\equiv( {\hat\Lambda- \Lambda^* })' \hat\Sigma_{ee}^{ - 1}
\hat\Lambda\hat H$ is defined following~(\ref{pthm10}). The first
term is $\|N^{1/2}\hat H^{1/2}\|\cdot O_p(T^{-1/2})$ by Lemma A.3(b)~\cite{bai2011Supplement}
and the second term is
\[
A\cdot\|N^{1/2}\hat H^{1/2}\|\cdot
O_p(T^{-1/2}).
\]
The fourth term is the transpose of the fifth. The
sixth term is given in Lem\-ma~A.3(d). The seventh term is bounded by $\|
\hat H\|\cdot\|{\sum_{i=1}^N}\frac{1}{\hat\sigma_i^4}\hat\lambda
_i\hat
\lambda_i'(\hat\sigma_i^2-\sigma_i^2)\hat H\|$. The term $\|{\sum
_{i=1}^N}\frac{1}{\hat\sigma_i^4}\hat\lambda_i\hat\lambda_i'(\hat
\sigma
_i^2-\sigma_i^2)\hat H\|$ is bounded by
$2C^4\sqrt r$ due to $|\frac{1}{\hat\sigma_i^2}(\hat\sigma
_i^2-\sigma
_i^{*2})|\le2C^4$ because of the boundedness of $\hat\sigma
_i^2,\sigma
_i^{*2}$. So the seventh term is $o_p(1)$ by~(\ref{pthm9}).
Given these results, in terms of $A$, equation~(\ref{pthm21}) can be
rewritten as
%
%
\begin{eqnarray}\label{august}\quad
\hat M_{f f}-M_{f f}^*&=&-A'M_{f
f}^*-M_{f f}^*A+A'M_{f f}^*A+\|N^{1/2}\hat H^{1/2}\|\cdot O_p(T^{-1/2})
\nonumber\\
&&{}-A\cdot\|N^{1/2}\hat H^{1/2}\|\cdot O_p(T^{-1/2})\\
&&{}+\|N^{1/2}\hat
H^{1/2}\|^2\cdot O_p(T^{-1/2})+o_p(1).\nonumber
\end{eqnarray}
However, by the definition of $\hat H$, $N\hat H=(\frac{1}{N}\hat
\Lambda
'\hat\Sigma_{ee}^{-1}\hat\Lambda)^{-1}$. Equation~(\ref{pthm17}) yields
$(\frac{1}{N}\hat\Lambda'\hat\Sigma_{ee}^{-1}\hat\Lambda
)^{-1}=(I_r-A)'(\frac{1}{N}\Lambda^{*\prime}\Sigma_{ee}^{*-1}\Lambda
^*)^{-1}(I_r-A)+o_p(\|I_r-A\|^2)$.
So we have
\begin{eqnarray*}
\|N^{1/2}\hat H^{1/2}\|^2&=&\tr[N\hat H]\\
&=&\tr\biggl[(I_r-A)'\biggl(\frac
{1}{N}\Lambda^{*\prime}\Sigma_{ee}^{*-1}\Lambda^*\biggr)^{-1}(I_r-A)+o_p(\|
I_r-A\|
^2)\biggr].
\end{eqnarray*}
The right-hand side is at most $O_p(A^2)$, implying that $\|N^{1/2}\hat
H^{1/2}\|=O_p(A)$.

Given the above result, we argue that the matrix $A$ must be
stochastically bounded. First, notice that the left-hand side of~(\ref{august}) is stochastically bounded by Assumption~\ref{AssumptionD}.
So if\vspace*{1pt} $A$ is not
stochastically bounded, the right-hand side is dominated by $A'M_{f
f}^*A$ in view of $\|N^{1/2}\hat H^{1/2}\|=O_p(A)$, But $A'M_{f f}^*
A$ will be unbounded since $M_{f f}^*$ is positive definite. A~contradiction is obtained.
Thus $A=O_p(1)$; it follows that $\|N^{1/2}\hat H^{1/2}\|=O_p(A)=O_p(1)$.
Given this result, we have
%
%
\begin{equation}\label{pthm79}
\hat M_{f f}-M_{f f}^*=-A'M_{f
f}^*-M_{f f}^*A+A'M_{f f}^*A+o_p(1).
\end{equation}

Next consider~(\ref{pthm22}). The seventh to ninth terms of the
right-hand\break side of~(\ref{pthm22}) are all $o_p(1)$ by Lemma A.3
\cite{bai2011Supplement} and $\|
N^{1/2}\hat H^{1/2}\|\,{=}\,O_p(1)$.
The~last term is bounded by
$2 C^3 \|\hat M_{f f}^{-1}\|\,{\cdot}\,\|\frac{1}{\hat\sigma_j}\hat
\lambda
_j\hat H^{1/2}\|\,{\cdot}\,\|H^{1/2}\|$.\vspace*{-1pt} Since\break ${\sum_{j=1}^N}\|\frac
{1}{\hat
\sigma_j}\hat\lambda_j\hat H^{1/2}\|^2\,{=}\,r$, $\|\hat M_{f f}^{-1}\|
=O_p(1)$ by Assumption~\ref{AssumptionD} and $\hat H=o_p(1)$ by~(\ref{pthm9}), it
follows that the last term is $o_p(1)$. The third and fourth terms are
similar to the fourth and fifth terms in~(\ref{pthm21}), and hence are
$o_p(1)$ due to $A=O_p(1), \|\lambda_j^*\|\le C$ for all $j$ and $\hat
M_{f f}$ being bounded.
Using the same arguments for~(\ref{pthm79}), we have
%
%
\begin{equation}\label{pthm80}
\hat\lambda_j-\lambda_j^*=\hat M_{f
f}^{-1}M_{f f}^*A\lambda_j^*- \hat M_{f f}^{-1}A'M_{f f}^*A\lambda
_j^* +o_p(1).
\end{equation}

We next prove consistency by using the
identification conditions.

\textit{Under} IC1: Since the identification condition is $\Lambda
^*=[I_r,\Lambda_2^{*\prime}]'$ and $\hat\Lambda=[I_r, \hat\Lambda
_2']'$, the
first $ r\times r$ upper block of $\Lambda^*$ is the same as that of
$\hat\Lambda$, that is, $[\lambda_1^*,\lambda_2^*,\ldots,\lambda
_r^*]=[\hat\lambda_1,\hat\lambda_2,\ldots,\hat\lambda_r]=I_r$.
By~(\ref{pthm80}) with $\hat\lambda_j-\lambda_j^*=0$, $j=1,2,\ldots,r$,
\[
M_{f f}^* A- A'M_{f f}^* A\xrightarrow{p}0.
\]
We now attach a~subscript to matrix $A$ to signify which identification
condition is used.
For IC$k$ ($k=1,2,\ldots,5$), we use $A_k$ to denote the corresponding $A$.
So the above equation implies $M_{f f}^*A_1-A_1'M_{f
f}^*A_1\xrightarrow{p}0$. Taking transpose, we have $A_1'M_{f
f}^*-A_1'M_{f f}^*A_1\xrightarrow{p}0$. Thus $M_{f f}^*A_1-A_1'M_{f
f}^*\xrightarrow{p}0$. Post-multiplying $A_1$, we obtain $M_{f
f}^*A_1^2-A_1'M_{f f}^*A_1\xrightarrow{p}0$. But we also have $M_{f
f}^*A_1-A_1'M_{f f}^*A_1\xrightarrow{p}0$. Thus $M_{f f}^*
A_1^2-M_{f f}^* A_1\xrightarrow{p}0$.\vadjust{\goodbreak} Since $M_{f f}^*$ is positive
definite, we have $A_1(I_r-A_1)\xrightarrow{p}0$. Since we have proved
that $I_r-A_1$ converges in probability to a~nonsingular matrix, it
follows that $A_1\xrightarrow{p}0$.

From~(\ref{pthm12}) and $A_1\xrightarrow{p}0$, we obtain $\frac{1}{N}(
{\hat\Lambda- \Lambda^* } )'\hat\Sigma_{ee}^{ - 1} ( {\hat\Lambda
- \Lambda^* } )\xrightarrow{p}0$, which is equivalent to~(\ref{sec411}).
From~(\ref{pthm79}) and $A_1\xrightarrow{p}0$, we obtain $\hat M_{f
f}-M_{f f}^*\xrightarrow{p}0$,
which is~(\ref{sec413}). This proves Proposition~\ref{thm1} under IC1.

\textit{Under} IC2: From the identification condition $\frac
{1}{N}\hat
\Lambda'\hat\Sigma_{ee}^{ - 1} \hat\Lambda= \frac{1}{N}\Lambda
^{*\prime}\Sigma_{ee}^{*-1}\Lambda^* = I_r $, by adding and subtracting
terms, we have the identity
%
%
\begin{eqnarray}\label{iden}
&&\frac{1}{N}(\hat\Lambda- \Lambda^*)'\hat\Sigma_{ee}^{ - 1} \hat
\Lambda+ \frac{1}{N}\hat\Lambda'\hat\Sigma_{ee}^{ - 1} (\hat
\Lambda- \Lambda^*)\nonumber\\[-8pt]\\[-8pt]
&&\qquad = - \frac{1}{N}\Lambda^{*\prime}( {\hat\Sigma_{ee}^{ - 1} - \Sigma
_{ee}^{*- 1} } )\Lambda^* + \frac{1}{N}( {\hat\Lambda- \Lambda^* }
)'\hat\Sigma_{ee}^{ - 1} ( {\hat\Lambda- \Lambda^* }).
\nonumber
\end{eqnarray}
By~(\ref{convar}) and Lemma A.4~\cite{bai2011Supplement}, the term $\frac{1}{N}\Lambda^{*\prime} (
{\hat\Sigma_{ee}^{ - 1} - \Sigma_{ee}^{*- 1} } )\Lambda^*$ is
$o_p(1)$. Thus
\[
\frac{1}{N}( {\hat\Lambda- \Lambda^* } )' \hat\Sigma_{ee}^{ - 1}
\hat\Lambda+ \frac{1}{N}\hat\Lambda' \hat\Sigma_{ee}^{ - 1} (
{\hat\Lambda- \Lambda^* } )-\frac{1}{N}( {\hat\Lambda- \Lambda^*
} )' \hat\Sigma_{ee}^{ - 1} ( {\hat\Lambda- \Lambda^* }
)\xrightarrow{p}0.
\]
The above can be written in terms of matrix $A$ (i.e., $A_2$ under IC2),
\[
A_2+A_2'-\frac{1}{N}( {\hat\Lambda- \Lambda^* } )' \hat\Sigma
_{ee}^{ - 1} ( {\hat\Lambda- \Lambda^* } )\xrightarrow{p}0.
\]
With $\frac{1}{N}\hat\Lambda'\hat\Sigma_{ee}^{ - 1} \hat\Lambda
=I_r$,~(\ref{pthm12}) implies $A_2A_2'-\frac{1}{N}( {\hat\Lambda-
\Lambda^* } )' \hat\Sigma_{ee}^{ - 1} ( {\hat\Lambda- \Lambda^* }
)\xrightarrow{p}0$. These two results imply
$A_2+A_2'-A_2A_2'\xrightarrow{p}0$, which is equivalent to
\[
(A_2-I_r)(A_2-I_r)'-I_r\xrightarrow{p}0.
\]
However,~(\ref{pthm79}) is equivalent to
\[
\hat M_{f f}-(A_2-I_r)'M_{f f}^*(A_2-I_r)\xrightarrow{p}0.
\]
Under IC2, $M_{f f}^*$ is a~diagonal matrix with distinct elements.
Also, $\hat M_{f f}$ is a~diagonal matrix by restriction. Applying
Lemma A.1~\cite{bai2011Supplement} with $Q=A_2-I_r$, $V=M_{f f}^*$, and $D=\hat M_{f f}$, we
conclude that $Q$ and thus $A_2-I_r$ converge in probability to a~diagonal matrix with elements being either $-1$ or 1. Equivalently,
$A_2$ converges to a~diagonal matrix with diagonal elements being
either 0 or 2. By assuming that $\hat\Lambda$ and $\Lambda^*$ have the
same column signs, we rule out 2 as the diagonal element So $A_2=o_p(1)$.
The rest of the proof is identical to IC1, implying Proposition~\ref{thm1}
under~IC2.

\textit{Under} IC3: IC3 requires $\hat M_{f f}=M_{f f}^*=I_r$, and so
by~(\ref{pthm79}), $(A_3-I_r)(A_3-I_r)'-I_r\xrightarrow{p}0$.
From~(\ref{pthm17}),
\[
\frac{1}{N}\Lambda^{*\prime}\Sigma_{ee}^{*-1}\Lambda^*-(A_3-I_r)'
\biggl(\frac{1}{N}\hat\Lambda'\hat\Sigma_{ee}^{-1}\hat\Lambda
\biggr)(A_3-I_r)\xrightarrow{p}0.
\]
Under IC3, $\frac{1}{N}\Lambda^{*\prime}\Sigma_{ee}^{*-1}\Lambda^*$ is
diagonal with distinct elements, and $\frac{1}{N}\hat\Lambda'\hat
\Sigma
_{ee}^{-1}\hat\Lambda$ is also diagonal by estimation restriction. The
latter matrix has distinct diagonal elements\vadjust{\goodbreak} with probability 1. It
follows that $A_3-I_r$ converges in probability to a~diagonal matrix
with diagonal elements either 1 or $-1$ by Lemma~A.1~\cite{bai2011Supplement} applied with
$Q=(A_3-I_r)'$, $V=\frac{1}{N}\hat\Lambda'\hat\Sigma_{ee}^{-1}\hat
\Lambda$, and $D=\frac{1}{N}\Lambda^{*\prime}\Sigma_{ee}^{*-1}\Lambda^*$.
The remaining proof is identical to that of IC2 and hence omitted. So
we have proved Proposition~\ref{thm1} under IC3.

\textit{Under} IC4: By the identification condition, both $\Lambda_1^*$
and $\hat\Lambda_1$ are lower triangular matrices, where $\Lambda_1$ is
first $r\times r$ submatrix of $\Lambda$.
Consider the first $r$ equations of~(\ref{pthm80}),
\[
\hat\Lambda_1' - \Lambda_1^{*\prime} - \hat M_{f f}^{ - 1} (
{M_{f
f}^* A_4 - A'_4 M_{f f}^* A_4 } )\Lambda_1^{*\prime}\xrightarrow{p}0.
\]
By~(\ref{pthm79}), we have
$\hat M_{f f} - M_{f f}^* + A_4' M_{f f}^* + M_{f f}^* A_4 - A'_4
M_{f f}^* A_4\xrightarrow{p}0$, which can be rewritten as
\[
\hat M_{f f}-(I_r-A_4)'M_{f f}^*(I_r-A_4)\xrightarrow{p}0.
\]
The above two equations imply
%
%
\begin{equation}\label{pthm33}\qquad
( {I_r - A~_4 })' M_{f f}^* ( {I_r - A~_4 } )( {\hat\Lambda_1' -
\Lambda_1^{*\prime} }) - ( {I_r - A~_4 } )' M_{f f}^* A_4 \Lambda
_1^{*\prime}\xrightarrow{p}0.
\end{equation}
Since both $\hat M_{f f}$ and $M_{f f}$ are of full rank, $\hat M_{f
f}-(I_r-A_4)'M_{f f}^*(I_r-A_4)\xrightarrow{p}0$ implies $I_r-A_4$ is
of full rank. Pre-multiplying $[(I_r-A_4)'M_{f f}^*]^{-1}$,
we obtain
%
%
\begin{equation}\label{pthm34}
(I_r-A_4)\hat\Lambda_1'-\Lambda_1^{*\prime}\xrightarrow{p}0.
\end{equation}
Since both $\hat\Lambda_1$ and $\Lambda_1^*$ are lower triangular with
diagonal elements all 1, both matrices are invertible. It follows that
$I_r-A_4-\Lambda_1^{*\prime}(\hat\Lambda_1')^{-1}\xrightarrow{p}0$. Since
both $\hat\Lambda_1$ and $\Lambda_1^*$ are lower triangular, we have
$I_r-A_4$ converges to an upper triangular matrix. However, $\hat M_{f
f}$ and $M_{f f}^*$ are both diagonal matrices and invertible. For
$\hat M_{f f}-(I_r-A_4)'M_{f f}^*(I_r-A_4)\xrightarrow{p}0$ to hold,
given that $I_r-A_4$ is an upper triangular matrix, it implies that
$I_r-A_4$ converges to a~diagonal matrix. Because both $\hat\Lambda_1$
and $\Lambda_1^*$ are matrices with diagonal elements 1, and given the
asymptotic diagonality of $A_4$, it follows by~(\ref{pthm34}) that
$I_r-A_4\xrightarrow{p}I_r$. So we have $A_4\xrightarrow{p}0$. The
remaining proof is the same as in IC1 and is omitted. This completes
the proof for IC4.

\textit{Under} IC5: 
Both $\hat M_{f f}$ and $M_{f f}^*$
are identity matrices; it follows from~(\ref{pthm79}) that
$(I_r-A_5)'(I_r-A_5)-I_r\xrightarrow{p}0$. The derivation of~(\ref{pthm34})
only involves the full rank of $I_r-A$, so it is applicable for IC5,
that is, $(I_r-A_5)\hat\Lambda_1'-\Lambda_1'\xrightarrow{p}0$. Since
both $\hat\Lambda_1$ and $\Lambda_1$ are lower triangular and
invertible, it follows that $I_r-A_5$ converges to an upper triangular
matrix. Given this result and $(I_r-A_5)'(I_r-A_5)-I_r\xrightarrow
{p}0$, it follows that $A_5$ converges to a~diagonal matrix with
diagonal elements either 0 or~2. By assuming that the column signs of
$\hat\Lambda$ and $\Lambda^*$ are the same,
we have $A_5\xrightarrow{p}0$. The remaining proof is the same as in
IC1 and is omitted. So we have proved Proposition~\ref{thm1} under~IC5.
This completes the proof of Proposition~\ref{thm1}.\vadjust{\goodbreak}
\end{pf*}

The proofs for other results are provided in the supplement
\cite{bai2011Supplement}.
\end{appendix}

\section*{Acknowledgments}

The authors thank two anonymous referees, an Associate Editor and an
Editor for constructive comments.

\begin{supplement}[id=suppA]
\stitle{Supplement to ``Statistical analysis of
factor models of high dimension''}
\slink[doi]{10.1214/11-AOS966SUPP} 
\sdatatype{.pdf}
\sfilename{aos966\_supp.pdf}
\sdescription{In this supplement we provide the detailed proofs for
Theorems~\ref{thm2}--\ref{thm5} and~\ref{thm7}. We also give a~simple
and direct proof that the EM solutions satisfy the first order conditions.
Remarks are given on how to make use of matrix properties to write
a~faster computer program.}
\end{supplement}


\printaddresses


\begin{thebibliography}{28}


\bibitem{r2}
\begin{barticle}[mr]
\bauthor{\bsnm{Amemiya},~\bfnm{Yasuo}\binits{Y.}},
  \bauthor{\bsnm{Fuller},~\bfnm{Wayne~A.}\binits{W.~A.}} \AND
  \bauthor{\bsnm{Pantula},~\bfnm{Sastry~G.}\binits{S.~G.}}
(\byear{1987}).
\btitle{The asymptotic distributions of some estimators for a~factor analysis
  model}.
\bjournal{J. Multivariate Anal.}
\bvolume{22}
\bpages{51--64}.
\bid{doi={10.1016/0047-259X(87)90074-1}, issn={0047-259X}, mr={0890881}}
\bptok{imsref}%
\end{barticle}
\endbibitem

\bibitem{r3}
\begin{bbook}[mr]
\bauthor{\bsnm{Anderson},~\bfnm{T.~W.}\binits{T.~W.}}
(\byear{2003}).
\btitle{An Introduction to Multivariate Statistical Analysis}, \bedition{3rd}
  ed.
\bpublisher{Wiley}, \baddress{Hoboken, NJ}.
\bid{mr={1990662}}
\bptok{imsref}%
\end{bbook}
\endbibitem

\bibitem{r4}
\begin{barticle}[mr]
\bauthor{\bsnm{Anderson},~\bfnm{T.~W.}\binits{T.~W.}} \AND
  \bauthor{\bsnm{Amemiya},~\bfnm{Yasuo}\binits{Y.}}
(\byear{1988}).
\btitle{The asymptotic normal distribution of estimators in factor analysis
  under general conditions}.
\bjournal{Ann. Statist.}
\bvolume{16}
\bpages{759--771}.
\bid{doi={10.1214/aos/1176350834}, issn={0090-5364}, mr={0947576}}
\bptok{imsref}%
\end{barticle}
\endbibitem

\bibitem{r5}
\begin{binproceedings}[mr]
\bauthor{\bsnm{Anderson},~\bfnm{T.~W.}\binits{T.~W.}} \AND
  \bauthor{\bsnm{Rubin},~\bfnm{Herman}\binits{H.}}
(\byear{1956}).
\btitle{Statistical inference in factor analysis}.
In \bbooktitle{Proceedings of the {T}hird {B}erkeley {S}ymposium on
  {M}athematical {S}tatistics and {P}robability, 1954--1955, Vol. {V}}
\bpages{111--150}.
\bpublisher{Univ. California Press}, \baddress{Berkeley}.
\bid{mr={0084943}}
\bptok{imsref}%
\end{binproceedings}
\endbibitem

\bibitem{r6}
\begin{barticle}[mr]
\bauthor{\bsnm{Bai},~\bfnm{Jushan}\binits{J.}}
(\byear{2003}).
\btitle{Inferential theory for factor models of large dimensions}.
\bjournal{Econometrica}
\bvolume{71}
\bpages{135--171}.
\bid{doi={10.1111/1468-0262.00392}, issn={0012-9682}, mr={1956857}}
\bptok{imsref}%
\end{barticle}
\endbibitem

\bibitem{bai2011Supplement}
\begin{bmisc}[auto:STB|2012/02/22|15:12:06]
\bauthor{\bsnm{Bai},~\bfnm{J.}\binits{J.}} \AND
  \bauthor{\bsnm{Li},~\bfnm{K.}\binits{K.}}
(\byear{2012}).
\bhowpublished{Supplement to ``Statistical analysis of factor models of high
  dimension.'' DOI:\doiurl{10.1214/11-AOS966SUPP}.}
\bptok{imsref}%
\end{bmisc}
\endbibitem

\bibitem{r8}
\begin{barticle}[mr]
\bauthor{\bsnm{Bai},~\bfnm{Jushan}\binits{J.}} \AND
  \bauthor{\bsnm{Ng},~\bfnm{Serena}\binits{S.}}
(\byear{2002}).
\btitle{Determining the number of factors in approximate factor models}.
\bjournal{Econometrica}
\bvolume{70}
\bpages{191--221}.
\bid{doi={10.1111/1468-0262.00273}, issn={0012-9682}, mr={1926259}}
\bptok{imsref}%
\end{barticle}
\endbibitem

\bibitem{bai-ng-2010}
\begin{bmisc}[auto:STB|2012/02/22|15:12:06]
\bauthor{\bsnm{Bai},~\bfnm{J.}\binits{J.}} \AND
  \bauthor{\bsnm{Ng},~\bfnm{S.}\binits{S.}}
(\byear{2010}).
\bhowpublished{Principal components estimation and identification of the
  factors. Unpublished manuscript, Columbia Univ.}
\bptok{imsref}%
\end{bmisc}
\endbibitem

\bibitem{r9}
\begin{bmisc}[auto:STB|2012/02/22|15:12:06]
\bauthor{\bsnm{Breitung},~\bfnm{J.}\binits{J.}} \AND
  \bauthor{\bsnm{Tenhofen},~\bfnm{J.}\binits{J.}}
(\byear{2008}).
\bhowpublished{GLS estimation of dynamic factor models.
Working paper, Univ.
  Bonn.}
\bptok{imsref}%
\end{bmisc}
\endbibitem

\bibitem{r18}
\begin{bbook}[auto:STB|2012/02/22|15:12:06]
\bauthor{\bsnm{Campbell},~\bfnm{J.~Y.}\binits{J.~Y.}},
  \bauthor{\bsnm{Lo},~\bfnm{A.~W.}\binits{A.~W.}} \AND
  \bauthor{\bsnm{MacKinlay},~\bfnm{A.~C.}\binits{A.~C.}}
(\byear{1997}).
\btitle{The Econometrics of Financial Markets}.
\bpublisher{Princeton Univ. Press}, \baddress{Princeton, NJ}.
\bptok{imsref}%
\end{bbook}
\endbibitem

\bibitem{r10}
\begin{barticle}[mr]
\bauthor{\bsnm{Chamberlain},~\bfnm{Gary}\binits{G.}} \AND
  \bauthor{\bsnm{Rothschild},~\bfnm{Michael}\binits{M.}}
(\byear{1983}).
\btitle{Arbitrage, factor structure, and mean-variance analysis on large asset
  markets}.
\bjournal{Econometrica}
\bvolume{51}
\bpages{1281--1304}.
\bid{doi={10.2307/1912275}, issn={0012-9682}, mr={0736050}}
\bptok{imsref}%
\end{barticle}
\endbibitem

\bibitem{r11}
\begin{bmisc}[auto:STB|2012/02/22|15:12:06]
\bauthor{\bsnm{Choi},~\bfnm{I.}\binits{I.}}
(\byear{2007}).
\bhowpublished{Efficient estimation of factor models. Working paper. Available
  at \texttt{
\href{http://hompi.sogang.ac.kr/inchoi/workingpaper/efficient\_estimation\_in\_choi\_et1.pdf}{http://hompi.sogang.ac.kr/inchoi/workingpaper/efficient\_estimation\_in\_}
\href{http://hompi.sogang.ac.kr/inchoi/workingpaper/efficient\_estimation\_in\_choi\_et1.pdf}{choi\_et1.pdf}}.}
\bptok{imsref}%
\end{bmisc}
\endbibitem

\bibitem{r12}
\begin{barticle}[auto:STB|2012/02/22|15:12:06]
\bauthor{\bsnm{Connor},~\bfnm{G.}\binits{G.}} \AND
  \bauthor{\bsnm{Korajczyk},~\bfnm{R.~A.}\binits{R.~A.}}
(\byear{1988}).
\btitle{Risk and return in an equilibrium APT: Application of a~new test
  methodology}.
\bjournal{Journal of Financial Economics}
\bvolume{21}
\bpages{255--289}.
\bptok{imsref}%
\end{barticle}
\endbibitem

\bibitem{r13}
\begin{bmisc}[auto:STB|2012/02/22|15:12:06]
\bauthor{\bsnm{Doz},~\bfnm{C.}\binits{C.}},
  \bauthor{\bsnm{Giannone},~\bfnm{D.}\binits{D.}} \AND
  \bauthor{\bsnm{Reichlin},~\bfnm{L.}\binits{L.}}
(\byear{2006}).
\bhowpublished{A~quasi-maximum likelihood approach for large approximate
  dynamic factor models. Discussion Paper 5724, CEPR.}
\bptok{imsref}%
\end{bmisc}
\endbibitem

\bibitem{r14}
\begin{barticle}[auto:STB|2012/02/22|15:12:06]
\bauthor{\bsnm{Forni},~\bfnm{M.}\binits{M.}},
  \bauthor{\bsnm{Hallin},~\bfnm{M.}\binits{M.}},
  \bauthor{\bsnm{Lippi},~\bfnm{M.}\binits{M.}} \AND
  \bauthor{\bsnm{Reichlin},~\bfnm{L.}\binits{L.}}
(\byear{2000}).
\btitle{The generalized dynamic-factor model: Identification and estimation}.
\bjournal{Rev. Econom. Statist.}
\bvolume{82}
\bpages{540--554}.
\bptok{imsref}%
\end{barticle}
\endbibitem

\bibitem{gewekezhou96}
\begin{barticle}[auto:STB|2012/02/22|15:12:06]
\bauthor{\bsnm{Geweke},~\bfnm{J.}\binits{J.}} \AND
  \bauthor{\bsnm{Zhou},~\bfnm{G.}\binits{G.}}
(\byear{1996}).
\btitle{Measuring the price of the arbitrage pricing theory}.
\bjournal{The Review of Financial Studies}
\bvolume{9}
\bpages{557--587}.
\bptok{imsref}%
\end{barticle}
\endbibitem

\bibitem{goyal-et-al}
\begin{barticle}[auto:STB|2012/02/22|15:12:06]
\bauthor{\bsnm{Goyal},~\bfnm{A.}\binits{A.}},
  \bauthor{\bsnm{Perignon},~\bfnm{C.}\binits{C.}} \AND
  \bauthor{\bsnm{Villa},~\bfnm{C.}\binits{C.}}
(\byear{2008}).
\btitle{How common are common return factors across the NYSE and Nasdaq?}
\bjournal{Journal of Financial Economics}
\bvolume{90}
\bpages{252--271}.
\bptok{imsref}%
\end{barticle}
\endbibitem


\bibitem{r17}
\begin{bbook}[mr]
\bauthor{\bsnm{Lawley},~\bfnm{D.~N.}\binits{D.~N.}} \AND
  \bauthor{\bsnm{Maxwell},~\bfnm{A.~E.}\binits{A.~E.}}
(\byear{1971}).
\btitle{Factor Analysis as a~Statistical Method}, \bedition{2nd} ed.
\bpublisher{Elsevier}, \baddress{New York}.
\bid{mr={0343471}}
\bptok{imsref}%
\end{bbook}
\endbibitem

\bibitem{NeMcF94}
\begin{bincollection}[mr]
\bauthor{\bsnm{Newey},~\bfnm{Whitney~K.}\binits{W.~K.}} \AND
  \bauthor{\bsnm{McFadden},~\bfnm{Daniel}\binits{D.}}
(\byear{1994}).
\btitle{Large sample estimation and hypothesis testing}.
In \bbooktitle{Handbook of Econometrics, {V}ol. {IV}}
(\beditor{R. F. Engle} and \beditor{D. McFadden}, eds.).
\bseries{Handbooks in Economics}
\bvolume{2}
\bpages{2111--2245}.
\bpublisher{North-Holland}, \baddress{Amsterdam}.
\bid{mr={1315971}}
\bptok{imsref}%
\end{bincollection}
\endbibitem

\bibitem{r19}
\begin{barticle}[mr]
\bauthor{\bsnm{Ross},~\bfnm{Stephen~A.}\binits{S.~A.}}
(\byear{1976}).
\btitle{The arbitrage theory of capital asset pricing}.
\bjournal{J. Econom. Theory}
\bvolume{13}
\bpages{341--360}.
\bid{issn={0022-0531}, mr={0429063}}
\bptok{imsref}%
\end{barticle}
\endbibitem

\bibitem{rubin1982algorithms}
\begin{barticle}[mr]
\bauthor{\bsnm{Rubin},~\bfnm{Donald~B.}\binits{D.~B.}} \AND
  \bauthor{\bsnm{Thayer},~\bfnm{Dorothy~T.}\binits{D.~T.}}
(\byear{1982}).
\btitle{E{M} algorithms for {ML} factor analysis}.
\bjournal{Psychometrika}
\bvolume{47}
\bpages{69--76}.
\bid{doi={10.1007/BF02293851}, issn={0033-3123}, mr={0668505}}
\bptok{imsref}%
\end{barticle}
\endbibitem

\bibitem{r20}
\begin{barticle}[mr]
\bauthor{\bsnm{Stock},~\bfnm{James~H.}\binits{J.~H.}} \AND
  \bauthor{\bsnm{Watson},~\bfnm{Mark~W.}\binits{M.~W.}}
(\byear{2002}).
\btitle{Forecasting using principal components from a~large number of
  predictors}.
\bjournal{J. Amer. Statist. Assoc.}
\bvolume{97}
\bpages{1167--1179}.
\bid{doi={10.1198/016214502388618960}, issn={0162-1459}, mr={1951271}}
\bptok{imsref}%
\end{barticle}
\endbibitem

\bibitem{r21}
\begin{barticle}[mr]
\bauthor{\bsnm{Stock},~\bfnm{James~H.}\binits{J.~H.}} \AND
  \bauthor{\bsnm{Watson},~\bfnm{Mark~W.}\binits{M.~W.}}
(\byear{2002}).
\btitle{Macroeconomic forecasting using diffusion indexes}.
\bjournal{J. Bus. Econom. Statist.}
\bvolume{20}
\bpages{147--162}.
\bid{doi={10.1198/073500102317351921}, issn={0735-0015}, mr={1963257}}
\bptok{imsref}%
\end{barticle}
\endbibitem

\bibitem{Vaart98}
\begin{bbook}[mr]
\bauthor{\bparticle{van~der} \bsnm{Vaart},~\bfnm{A.~W.}\binits{A.~W.}}
(\byear{1998}).
\btitle{Asymptotic Statistics}.
\bseries{Cambridge Series in Statistical and Probabilistic Mathematics}
\bvolume{3}.
\bpublisher{Cambridge Univ. Press}, \baddress{Cambridge}.
\bid{mr={1652247}}
\bptok{imsref}%
\end{bbook}
\endbibitem

\bibitem{wu1983convergence}
\begin{barticle}[mr]
\bauthor{\bsnm{Wu},~\bfnm{C.~F.~Jeff}\binits{C.~F.~J.}}
(\byear{1983}).
\btitle{On the convergence properties of the {EM} algorithm}.
\bjournal{Ann. Statist.}
\bvolume{11}
\bpages{95--103}.
\bid{doi={10.1214/aos/1176346060}, issn={0090-5364}, mr={0684867}}
\bptok{imsref}%
\end{barticle}
\endbibitem

\end{thebibliography}
\end{document}